\documentclass{article}
\usepackage[utf8]{inputenc}
\usepackage[utf8]{inputenc}
\usepackage{amsmath}
\usepackage[left=2.9cm, right=2.9cm, bottom=2.5cm, top=2.5cm]{geometry}
\usepackage{hyperref}
\usepackage{cleveref}
\usepackage{comment}
\usepackage{amssymb}
\usepackage{mathtools}
\usepackage{tikz-cd}
\usepackage{tikz}
\usetikzlibrary{arrows.meta}
\usetikzlibrary{bending}
\usetikzlibrary{decorations.text}
\usetikzlibrary{calc, arrows, decorations.markings}
\usepackage{amsthm}
\usepackage{graphicx}
\usepackage{graphics}
\usepackage{mathrsfs}
\usepackage{cleveref}
\usepackage{bigints}
\usepackage{thmtools}
\usepackage{enumitem}
\usepackage{changepage}
\usepackage{titlesec}
\usetikzlibrary{arrows}
\newcommand{\midarrow}{\tikz \draw[-triangle 90] (0,0) -- +(.1,0);}
\usetikzlibrary{arrows.meta}
\newlist{steps}{enumerate}{1}
\usepackage{amsmath,amssymb}
\setlist[steps, 1]{label = Step \arabic*:}
\theoremstyle{plain}
\usetikzlibrary{positioning, decorations.markings, decorations.pathreplacing}
\tikzstyle{directed}=[postaction={decorate,decoration={markings,mark=at position .5 with{\arrow{stealth}}}}]
\newtheorem{theorem}{Theorem}
\newtheorem{definition}{Definition}

\newtheorem{proposition}{Proposition}
\newtheorem{lemma}{Lemma}
\newtheorem{remark}{Remark}

\theoremstyle{definition}
\newtheorem{exmp}{Example}

\newtheorem{conjecture}{Conjecture}
\numberwithin{equation}{subsection}
\usepackage{rotating}
\usepackage{comment}
\usepackage[nottoc,notlot,notlof]{tocbibind} 

\usepackage{bbm}
\title{A Formal Equivalence of Deformation Quantization and Geometric Quantization (of Higher Groupoids) and Non-Perturbative Sigma Models}
\author{Joshua Lackman\footnote{jlackman@math.toronto.edu}}
\date{}

\begin{document}

\maketitle
\begin{abstract}
\noindent Based on work done by Bonechi, Cattaneo, Felder and Zabzine on Poisson sigma models, we formally show that Kontsevich's star product can be obtained from the twisted convolution algebra of the geometric quantization of a Lie 2-groupoid, one which integrates the Poisson structure. We show that there is an analogue of the Poisson sigma model which is valued in Lie 1-groupoids and which can often be defined non-perturbatively; it can be obtained by symplectic reduction using the quantization of the Lie 2-groupoid. We call these groupoid-valued sigma models and we argue that, when they exist, they can be used to compute correlation functions of gauge invariant observables. This leads to a (possibly non-associative) product on the underlying space of functions on the Poisson manifold, and in several examples we show that we recover strict deformation quantizations, in the sense of Rieffel. Even in the cases when our construction leads to a non-associative product we still obtain a $C^*$-algebra and a ``quantization" map. In particular, we construct noncommutative $C^*$-algebras equipped with an $SU(2)$-action, together with an equivariant ``quantization" map from $C^{\infty}(S^2)\,.$ No polarizations are used in the construction of these algebras.
\end{abstract}
\tableofcontents
\section{Introduction}
\subsection{Deformation Quantization and the Poisson Sigma Model}$\,$
\vspace{0.05cm}\\
In his seminal paper \cite{kontsevich}, Kontsevich showed that every Poisson manifold $(M,\Pi)$ admits a (formal) deformation quantization, ie. there is an associative product (called a star product) $\star$ on $C^{\infty}(M)[[\hbar]]$\footnote{The ``star" in star product refers to the fact that the associative product is induced by a sequence of bidifferential operators, see \cite{kontsevich}.} which, formally, deforms the product on $C^{\infty}(M)$ in such a way that 
\vspace{0.05cm}\\
\begin{equation*}
    f*g-g*f=i\hbar\{f,g\}+\mathcal{O}(\hbar^2)\,.
    \end{equation*}
    \vspace{0.05cm}\\
In \cite{catt}, Cattaneo and Felder provide a Poisson sigma model approach to defining the star product. The Poisson sigma model is, in particular, a topological bosonic open string theory, and the authors show that the star product can be obtained by computing a 3-point correlation function of boundary fields.
\vspace{0.5cm}\\More precisely, the target space of the Poisson sigma model is a Poisson manifold $(M,\Pi)$ and there are two bosonic fields, given by a vector bundle morphism 
\vspace{0.05cm}\\
\begin{equation*}
    (X,\eta):TD\to T^*M
\end{equation*}
\vspace{0.05cm}\\
covering $X:D\to M\,,$ where $D$ is a disk with three cyclically ordered marked points on the boundary, denoted $0\,,1\,,\infty\,.$ The first fields are given by $X$ and the second fields are given by the $X^*(T^*M)$-valued one forms $\eta\,,$ subject to the boundary condition that the pullback of $\eta$ to $\partial D$ takes values in the zero section. The star product is then given by
\vspace{0.05cm}\\
\begin{equation}\label{star}
 (f\star g)(m)=\int_{X(\infty)=m} f(X(1))g(X(0))\,e^{\frac{i}{h}S[X,\eta]}\,dX\,d\eta\,.
\end{equation}
\vspace{0.05cm}\\
Here $m\in M\,,f\,,g\in C^{\infty}(M)\,,$ the integral is over all fields $X\,,\eta$ with $X(\infty)=m\,,$ and the action is given by
\vspace{0.05cm}\\
\begin{equation}\label{action}
    S[X,\eta]=\int_D \eta\wedge dX+\frac{1}{2}\Pi\vert_X(\eta\,,\eta)\,.
\end{equation}
\vspace{0.05cm}\\
The aforementioned authors show that the semiclassial expansion of this product about the classical solution $(X,\eta)=(m,0_m)$ is equivalent to Kontsevich's star product.
\subsection{An Equivalent Path Integral Over Lie Algebroid Morphisms}$\,$\\
Later, in \cite{bon}, after using polar coordinates on the disk and formally integrating out $\eta_{\theta}\,,$ Bonechi, Cattaneo and Zabzine (formally) show that \ref{star} is equivalent to a path integral with a different domain of integration, one which emphasizes the Lie algebroid $T^*M$ associated to $(M,\Pi)\,.$ Rather than considering all such vector bundle morphisms $(X,\eta)\,,$ we may consider only those morphisms which are morphisms of Lie algebroids. Now, let $X\in\text{Hom}_{\text{LA}}(TD, T^*M)$ denote Lie algebroid morphisms. The product \ref{star} can equivalently be written as 
\vspace{0.05cm}\\
\begin{align}\label{mystarr}
    (f\star g)(m)=\int_{\pi(X(\infty))=m}f(\pi(X(1)))\,g(\pi(X(0)))\,e^{\frac{i}{\hbar}S[X]}\,dX\,,
\end{align}
where $\pi:T^*M\to M$ is the projection and
\begin{equation}
S[X]=\int_D X^*\Pi\,.
\end{equation}
\vspace{0.05cm}\\
One nice feature of \ref{mystarr} is that the path integral defining $(f\star g)(m)$ only involves values of $f$ and $g$ on the symplectic leaf containing $m\,,$ in contrast to \ref{star}. This means that, formally, \ref{mystarr} is tangential, which essentially means that it restricts to symplectic leaves. For star products as in \cite{kontsevich}, this essentially means that each bidifferential operator is tangent to the leaves (see \cite{gamme}).\\
\subsection{The Main Idea}\label{idea}
In this paper we discuss a possible quantization scheme for Poisson manifolds, one which doesn't use polarizations in an essential way. We do this by taking a path integral perspective. Here we will summarize the idea:\\\\
\begin{enumerate}
    \item 
Let $(M,\Pi)$ be a Poisson manifold. Consider the simplicial space $V^{\bullet}$ which in degree $n$ is given by
\begin{equation*}
    V^n :=\text{Hom}_\text{VB}(\,T|\Delta^n|,\,T^*M),
\end{equation*}
ie. vector bundle morphisms $T|\Delta^n|\to T^*M,$ where $|\Delta^n|$ is the standard n-simplex (the space, not the simplicial set). We have that $V^1\cong P(T^*M),$ where for a space $A,$ $PA$ denotes the space of paths $[0,1]\to A.$ We have that $P(T^*M)\cong T^*(PM),$ therefore $V^1$ has a canonical symplectic form. So far we haven't used $\Pi.$
\item Associated with $\Pi$ is a Hamiltonian action on $V^1,$ and if $\mu$ is the moment map, then $\mu^{-1}(0)$ is the coisotropic submanifold consisting of Lie algebroid paths (where the Lie algebroid is the one associated with $\Pi,$ where the underlying vector bundle is $T^*M$). See \cite{felder}.
\item The space of Lie algebroid paths sits inside a simplicial space, denoted $\Pi_{\infty}(T^*M),$ which in degree $n$ is given by\\
\begin{equation*}
  \Pi_{\infty}^{(n)}(T^*M)=\{f\in \text{Hom}(\,T|\Delta^n|,\, T^*M\,): f\vert_{\text{vertices}}=0\}\,,
\end{equation*}$\,$\\
See \cite{zhuc}.
\item The reduced phase space of $V^1$ is the space of arrows of the source simply connected symplectic groupoid integrating $\Pi,$ denoted $\Pi_1(T^*M).$ In other words, we can think of $\Pi_1(T^*M)$ as being the symplectic reduction of $V^{\bullet}.$ See \cite{felder}
\item Now since $V^1$ comes with a canonical symplectic form $\omega$ we can geometrically quantize $V^{\bullet}$ (ie. we can construct a multiplicative line bundle with connection prequantizing the $\omega$). When we pullback this geometric quantization to $\Pi_{\infty}(T^M)$ we get a ``geometric quantization" of $\Pi_{\infty}(T^*M),$ where the presymplectic form (ie. closed 2-form) is obtained by pulling back $\omega.$
\item After doing this we can form the twisted convolution algebra of $\Pi_{\infty}(T^*M),$ which for functions $s_1, s_2\in C^{\infty}(\Pi_{\infty}^{(1)}(T^*M))$ will be denoted $s_1\ast s_2,$ and is given by
\begin{equation}
    (s_1\ast s_2)(\gamma)=\int_{\begin{subarray}{l}X\in \Pi^{(2)}_{\infty}(T^*M)\\d_1(X)=\gamma\end{subarray}}s_1(d_0(X))s_2(d_2(X))e^{\frac{i}{\hbar}S[X]}\,dX\,,
\end{equation}
where $d_0, d_1, d_2$ are the face maps $\Pi_{\infty}^{(2)}(T^*M)\to \Pi_{\infty}^{(1)}(T^*M).$
\item There is a canonical ``quantization" map 
\begin{align*}
& C^{\infty}(M)\to C^{\infty}(\Pi_{\infty}^{(1)}(T^*M))\,,
 \\&f\mapsto \hat{f}\,,\;\hat{f}(\gamma)=f(\pi(\gamma(1/2)))\,,
\end{align*} 
where $\pi:T^*M\to M$ is the projection. It follows that $(\hat{f}\ast\hat{g})(m)$ agrees with \ref{mystarr}, and therefore this essentially gives the Poisson sigma model (formally).
\item Now there is a Lie 2-groupoid present, it is given by a quotient of $\Pi_{\infty}(T^*M)$ and will be denoted $\Pi_{2}(T^*M),$ see \cite{zhuc}. The objects and arrows of $\Pi_{2}(T^*M)$ are the same as those in $\Pi_{\infty}(T^*M),$ and we can perform a partial reduction to get a ``geometric quantization" of $\Pi_{2}(T^*M).$ The twisted convolution algebra of $\Pi_{2}(T^*M)$ therefore also agrees with \ref{mystarr}.
\item Using $\Pi_{2}(T^*M),$ we can think of the Poisson sigma model as a higher sigma model, valued in a Lie 2-groupoid. In order to do this we need to use the following identification, given by the Yoneda lemma:\\
\begin{equation*}
    \Pi_2^{(2)}(T^*M)\cong \text{Hom}(\Delta^2, \Pi_2(T^*M))\,.
\end{equation*}$\,$\\
Here $\Delta^2$ is the standard 2-simplex (this time it is the simplicial set). We think of this as a sigma model in the following way: the source and target are $\Delta^2,\, \Pi_2(T^*M),$ respectively, and the action is given by the one in \ref{mystarr}. Then we can write \ref{mystarr} in the following way: \\\begin{equation}\label{bah}
            (f\star g)(m)=\int_{\begin{subarray}{l}F\in\text{Hom}(\Delta^2,\Pi_2(T^*M)),\\F(1,1)=m\end{subarray}} dF_m\,e^{\frac{i}{\hbar}S(F(2,1),F(0,1))}\,f(F(2,2))g(F(0,0)) \,.
              \end{equation}$\,$\\
Here $dF_m$ is a measure on the space $\{F\in\text{Hom}(\Delta^2,\Pi_2(T^*M)): F(1,1)=m\}\,.$ We think of $m$ as being the constant algebroid path mapping to $m\,.$    
\item Now assume that \\
\begin{equation}
    \int_{S^2} X^*\Pi=0,
\end{equation} $\,$\\
for all Lie algebroid morphisms $X:TS^2\to T^*M.$ Then under symplectic reduction, the geometric quantization of $V^{\bullet}$ can be reduced to a geometric quantization of $\Pi_1(T^*M),$ and the higher sigma model can be reduced to a sigma model valued in $\Pi_1(T^*M),$ where the product is analogous to \ref{bah}. It is given by  \\
\begin{equation}\label{co}
            (f\star g)(m)=\int_{\begin{subarray}{l}F\in\text{Hom}(\Delta^2,\Pi_1(T^*M)),\\F(1,1)=m\end{subarray}} dF_m\,e^{\frac{i}{\hbar}S(F(2,1),F(0,1))}\,f(F(2,2))g(F(0,0)) \,.
              \end{equation}$\,$\\
The source and target spaces are $\Delta^2,\,\Pi_1(T^*M),$ respectively; $S$ is the analogue of the action, it is given by a $2$-cocycle on $\Pi_1(T^*M)$ which maps to $\Pi$ under the van Est map, see \ref{appen} (recall that, $\Pi$ itself defines a 2-cocycle on the Lie algebroid $T^*M$). Note that, 
\begin{equation*}
    \text{Hom}(\Delta^2,\Pi_1(T^*M))\cong \Pi_1^{(2)}(T^*M),
    \end{equation*}
    and therefore we are integrating over a finite dimensional space.
\item We can then use the natural volume form on each symplectic leaf to obtain a $C^*$-algebra, by using the operator norm of $C^{\infty}_c(M)$ acting on itself via \ref{co}. In particular, this $C^*$-algebra comes with a natural representation, which we think of as being our space of states (although not all of these states are pure, ie. they don't all live in irreducible representations). However, before doing this we need to take care of any possible divergences, which in our examples arise just due to the isotropy of $\Pi.$ Since we are integrating over a finite dimensional space, this can often be done non-perturbatively, as in the case of the Moyal product, which we do in section \ref{fund}.
\item We do not claim that this product is associative, however we will see several examples in which it is. However, even when the product isn't associative we still get a $C^*$-algebra, and thus it may still give a strict quantization of $\Pi,$ in the sense of Rieffel, see \cite{rieffel3}. Relevant to this is the fact that Kontsevich's proof of associativity of the star product in \cite{kontsevich} only shows that \ref{star} is \textit{perturbatively associative} when expanding about the constant critical point.\footnote{We are not claiming that \ref{star} is non-perturbatively non-associative.}
\item One point to note: the products \ref{co}, \ref{mystarr} are manifestly tangential, meaning that they pulback to give ``quantizations" of each symplectic leaves.
\item Graphically, the product \ref{co} (and \ref{mystarr}) can be represented in the following way: \[
\begin{tikzpicture}
\tikzstyle arrowstyle=[scale=2]  \tikzstyle directed=[postaction={decorate,decoration={markings,mark=at position .5 with{\arrow[arrowstyle]{stealth}}}}]
\foreach \i in {0, 180, 270}{
\filldraw (\i:1.25)circle(0.05);
}
\draw[dashed, directed] (0:1.25)arc[start angle=0, end angle=180, x radius=1.25cm, y radius=1.25cm] node[midway, above] {$\gamma_2\gamma_1\;\;\,$};
\draw[directed] (270:1.25)arc[start angle=-90, end angle=-180, x radius=1.25cm, y radius=1.25cm] node[midway, left] {$\gamma_2$};
\draw[directed] (0:1.25)arc[start angle=0, end angle=-90, x radius=1.25cm, y radius=1.25cm] node[midway, right] {$\;\,\gamma_1$};

      \path
      (-2,0) node[] {$\displaystyle\bigint$}
      (-3.5, 0) node[] {$(f\star g)(m)\;=$}
      (-1.25, 0) node[left] {$f$}
      (1.25, 0) node[right] {$g$}
      (0, -1.25) node[below] {$m$}
      (0, 0) node[] {$e^{\frac{i}{\hbar}S[\gamma_2,\gamma_1]}$}
      (-2.25, -1.2) node[] {$\gamma_1,\gamma_2$}
      
      ;
\end{tikzpicture}
\]$\,$\\
Here $S$ is the action and $f, g$ are being evaluated at the target of $\gamma_2$ and the source of $\gamma_1,$ respectively. Therefore, in the spirit of Feynman these prodcuts are given by a sum over all paths, however the paths are ones which lie in a (higher) groupoid.
\end{enumerate}$\,$\\
Now in the rest of the paper we will give progressively more detail about this idea. Let us remark that the simplicial space $V^{\bullet}$ is sitting in the background, however we discussed it more for conceptual reasons and we won't be using it much.\\
\subsection{Definitions: Tangential Quantizations, Gauge Equivalences and $\mathbf{C^*}$-Algebras}\label{explain}
We are mostly interested in strict quantizations/deformation quantizations (ie. not formal products) and will now give some definitions:\footnote{We are not quite giving the most general definitions here, but we are giving what we will use.}
\vspace{0.05cm}\\
\begin{definition}\label{strict}(see \cite{rieffel})
Let $(M,\Pi)$ be a Poisson manifold. A strict deformation quantization of $(M,\Pi)$ is, for each $\hbar\in\mathbb{R}\,,$ a $C^*$-algebra containing $C^{\infty}_c(M)$ as a dense $^{*_{\hbar}}$-subalgebra, denoted 
\vspace{0.05cm}\\
\begin{equation*}
    (C^{\infty}_c(M),\,\star_{\hbar}, \,^{*_{\hbar}},\,||\cdot ||_{\hbar})\,,
    \end{equation*}
such that
\begin{enumerate}
    \item For each $f\in C^{\infty}_c(M)\,,$ the map $\hbar\to ||f||_{\hbar}$ is continuous.
    \item For $\hbar=0$ the $C^*$-algebra is the natural one, ie. $(C^{\infty}_c(M),\cdot, ^*,||\cdot ||_{L^{\infty}(M)})\,,$
    where the product is multiplication and $^*$ denotes complex conjugation.
    \item $||(f\star_{\hbar}g-g\star_{\hbar}f)/\hbar-i\{f,g\}||_{\hbar}\xrightarrow[]{\hbar\to 0}0\,.$
\end{enumerate}
\end{definition}$\,$
\vspace{0.05cm}\\
Strict deformation quantizations are ideal, but a weaker notion, called strict quantization, is also interesting. The main difference is that a strict quantization assigns to each $f\in C^{\infty}_c(M)$ an element of a $C^*$-algebra, but its image doesn't need to be closed under the product in the $C^*$-algebra:
\begin{definition}\label{strict quant}(see \cite{rieffel3})
Let $(M,\Pi)$ be a Poisson manifold. A strict quantization of $(M,\Pi)$ is, for each $\hbar\in\mathbb{R}\,,$ a $C^*$-algebra $\mathcal{A}_{\hbar}$ together with a linear map $\mathcal{T}_{\hbar}:C^{\infty}_c(M)\to \mathcal{A}_{\hbar}$ (called a quantization map, which usually preserves the involution) such that $\mathcal{A}_{\hbar}$ is generated by the image of $C^{\infty}_c(M)\,,$ and such that
\begin{enumerate}
    \item For each $f\in C^{\infty}_c(M)\,,$ the map $\hbar\to ||\mathcal{T}_{\hbar}(f)||_{\hbar}$ is continuous.
    \item $\mathcal{A}_0=C^{\infty}_c(M)$ with its natural $C^*$-algebra structure and $\mathcal{T}_0$ is the identity.
    \item $||(\mathcal{T}_{\hbar}(f)\star_{\hbar}\mathcal{T}_{\hbar}(g)-\mathcal{T}_{\hbar}(g)\star_{\hbar}\mathcal{T}_{\hbar}(f))/\hbar-iT_{\hbar}(\{f,g\})||_{\hbar}\xrightarrow[]{\hbar\to 0}0\,.$
\end{enumerate}
\end{definition}$\,$
\vspace{0.05cm}\\
Now we will define what it means for a quantization to be tangential:
\vspace{0.05cm}\\
\begin{definition}\label{tang}(see \cite{rieffel3}
A strict deformation quantization of a Poisson manifold $(M,\Pi)$ is tangential, if for each symplectic
leaf $\mathcal{L}\,,$ the ideal
\vspace{0.05cm}\\
\begin{equation*}
    I_\mathcal{L}=\{f\in C^{\infty}_c(M) : f\vert_{\mathcal{L}}=0\}
\end{equation*}
\vspace{0.05cm}\\
of $C^{\infty}_c(M)$ is also an ideal for the $C^*$-algebra, for all $\hbar\,.$ If $\mathcal{L}$ is not closed, we further require that the closure of $\mathcal{L}$ consists of all points where all of the functions in $I_{\mathcal{L}}$ vanish. 
\vspace{0.5cm}\\Similarly, we say a strict quanzitzation is tangential if the $C^*$-subalgebra generated by  $\mathcal{T}_{\hbar}(I_{\mathcal{L}})$ is an ideal, for all $\hbar\,.$
\end{definition}$\,$
\vspace{0.05cm}\\
Essentially, tangential (strict) quantizations allow us to quantize each symplectic leaf separately. Now classically, a particle on a Poisson manifold is confined to a symplectic leaf, therefore tangential (strict) quantizations are exactly those which don't exhibit quantum tunneling, in this sense. Related to tangential deformation quantizations is a conjecture of Weinstein:\\
\begin{conjecture}\label{weinc}(see in \cite{wein})
The dual of a Lie algebra, $\mathfrak{g}^*,$ admits a tangential star product if and only if there is a flat torsion-free affine connection on a neighborhood of the identity in $G$ for
which the induced connection on the vector bundle $T^*G$ is compatible with the
decomposition into the left (or right) translates of coadjoint orbits.
\end{conjecture}
We will discuss the relevance of this conjecture at the end of this paper, but it has implications about associativity. 
\paragraph{Gauge Equivalences}$\,$\\
\vspace{0.05cm}\\
Next, the authors of \cite{bon} define gauge equivalences in such a way that reduced phase space of the Poisson sigma model is the symplectic groupoid (if it exists). Furthermore, two algebroid morphims $X,\,X':TD\to T^*M$ are gauge equivalent if they agree on the boundary. Using this notion of gauge equivalence, in order to make sense of \ref{mystarr}, a necessary condition is that
\begin{equation}\label{cond}
    \int_{S^2} X^*\Pi\in 2\pi\hbar\mathbb{Z}\,,
\end{equation} for any algebroid morphism $X:TS^2\to T^*M\,.$ The authors go on to speculate that in this case a non-perturbative definition of \ref{mystarr} should exist, and if $\Pi_1(T^*M)\rightrightarrows M$ is the symplectic groupoid, then the path integral should be defined over the space of maps 
\begin{equation*}
    [0,1]\to \Pi_1^{(1)}(T^*M)\,,
\end{equation*}
where  $\Pi_1^{(1)}(T^*M)$ is the space of arrows of $\Pi_1(T^*M)\rightrightarrows M\,.$
\vspace{0.5cm}\\
What we do in this paper is different than what is speculated in the previous paragraph — we are essentially going to be doing a symplectic reduction of the Poisson sigma model. Rather than two morphisms  $X,\,X':TD\to T^*M$ being gauge equivalent if they agree on the boundary, we define them to be gauge equivalent if they agree on the boundary, up to Lie algebroid homotopy. For this to work, we require a slightly stronger condition than \ref{cond}: we require that for any algebroid morphism $X:TS^2\to T^*M\,,$
\begin{equation}\label{cond2}
    \int_{S^2} X^*\Pi=0
\end{equation} 
(we will speculate about the more general case at the end of the paper). Then, rather than integrating over the infinite dimensional space of maps $[0,1]\to \Pi_1^{(1)}(T^*M)\,,$ we integrate over the finite dimensional space of composable pairs of arrows in $\Pi_1(T^*M)\,.$ To motivate this, let us remark here that one can interpret the fields, ie. $X\in\text{Hom}_{\text{LA}}(TD, T^*M)\,,$ as composable pairs of arrows (ie. a 2-simplex) in a higher integration of $(M,\Pi)\,.$ We use this to compute correlation functions in the Poisson sigma model with respect to gauge invariant observables, non-perturbatively.
\paragraph{The $\mathbf{C^*}$-Algebra}$\,$\\
\vspace{0.05cm}\\Suppose we have a reasonable, non-perturbative definition of \ref{mystarr}. Then, taking advantage of its tangential nature, we get a family of $C^*$-algebras (one for each $\hbar)$ associated to each symplectic leaf $\mathcal{L}$ by using the operator norm $||\cdot||^{\mathcal{L}}_{\hbar}$ of the action of $L^2(\mathcal{L})$ on itself (where the action is given by \ref{mystarr}), with respect to the canonical volume form on each symplectic leaf, and where the involution is complex conjugation. Then, we get the $\hbar$-family of $C^*$-algebras we are interested in as follows: \\
\begin{definition}
For $f\in C^{\infty}_c(M)\,,$ let $||f||_{\hbar}:=\text{sup}_{\mathcal{L}}\mathbf{||}\,f\vert_{\mathcal{L}}\,||^{\mathcal{L}}_{\hbar}\,.$ The involution is given by complex conjugation and the product is given by the product of operators.
    \end{definition}$\,$\\
 Essentially, all we're doing is taken the direct sum of all of the Hilbert spaces (associated with each symplectic leaf) and taking the $C^*$-algebra to be the one associated to the operator norm. In particular, we get a representation of this $C^*$-algebra on the Hilbert space associated each symplectic leaf.
\vspace{0.5cm}\\Before getting into all of details, we will further explain our idea in the next section. In particular, there we will explain how our quantization prodecure recovers the usual strict deformation quantization of symplectic vector spaces, ie. the integral form of the Moyal product. Let us emphasize here that our derivation is different than other derivations in the literature, in particular no polarizations or Fourier transforms will be used.
\section{An Explanation of This Paper and Essential Concepts}In this paper we are going to show that \ref{mystarr} can be obtained by using the twisted convolution algebra of a higher integration of $(M,\Pi)$; it suffices to use the Lie 2-groupoid in \cite{zhuc}, which we will denote by $\Pi_2(T^*M)\,.$ To obtain the twist we are going to perform the analogue of geometric quantization on $\Pi_2(T^*M)\,,$ which involves constructing an $\mathbb{R}$-valued 2-cocycle on $\Pi_2(T^*M)\,.$ This is very much inline with the geometric quantization program developed by Weinstein, Xu, Hawkins, which also involves integrating a Poisson manifold and forming a twisted convolution algebra (on a reduced groupoid); see \cite{weinstein}, \cite{weinstein1}, \cite{eli}. See section \ref{idea} for a brief summary of this paper.
\vspace{0.5cm}\\There is one main difference between how geometric quantization is usually performed and how it will be done here: we will not be needing a polarization. This is a benefit because nonsingular polarizations are notoriously difficult to find in this context. As a result, the emphasis will be placed on the 2-cocycle and the van Est map, rather than on the 2-form which will take a back seat. We are using the term ``geometric quantization" a bit loosely because (aside from the absense of a polarization) the 2-form on $\Pi_2(T^*M)$ is presymplectic (ie. closed), not symplectic — however, it does sit inside a symplecic manifold as a coisotropic manifold.
\vspace{0.5cm}\\After doing this we will explain an analogue of \ref{mystarr} on $\Pi_1(T^*M)$ (the source simply connected groupoid), which we call a groupoid-valued sigma model, and in several examples we will show that this leads to a genuine associative product which gives a strict deformation quantization. The example we flesh out the most, and the example we think is the most instructive, is the one done in \ref{fund}. Basically, the relationship goes: under symplectic reduction,\\
\begin{align*}
    &\text{Geometric quantization of }\Pi_2(T^*M)\leadsto\text{Geometric quantization of }\Pi_1(T^*M)\,,\\
    &f\star g \text{ with respect to }\Pi_2(T^*M) \leadsto f\star g \text{ with respect to }\Pi_1(T^*M)\,,
\end{align*}
\vspace{0.25cm}\\The appendix covers some topics not covered in the main text, mostly things related to the cohomology of Lie groupoids, Lie algebroids and the van Est map.
\vspace{0.5cm}\\Many authors (aside from the ones mentioned throughout this paper) have worked on related topics, eg. convergent products, see  \cite{amar}, \cite{ale}, \cite{gamme}, \cite{heinz} \cite{rieffel2}, \cite{Waldmann}. Also, \cite{rieffel3} contains a list of open problems\footnote{Some of them may have been resolved by now.} regarding exactly the topics in this paper.
\subsubsection{Geometric Quantization of Poisson Manifolds Using Lie Groupoids}
\paragraph{Twisted Convolution Algebra}$\,$\\
\vspace{0.05cm}\\To explain a bit of the preceding paragraphs, let us first recall twisted convolution algebras on Lie groupoids (when we say Lie groupoids we will always mean Lie 1-groupoids). 
\vspace{0.5cm}\\Let $G\rightrightarrows M$ be a Lie groupoid. We let $G^{(0)}=M$ and we let $G^{(n)}$ be the space of composable $n$-tuples of arrows; together these form a simplicial set. There is a differential $\delta^*:C^{\infty}(G^{(n)})\to C^{\infty}(G^{(n+1)})$ given by 
\vspace{0.05cm}\\
\begin{equation}
    \delta^*f=\sum_{i=0}^n (-1)^id_i^*f\,,
\end{equation}
\vspace{0.05cm}\\
where $d_i:G^{(n+1)}\to G^{(n)},\,i=0,\ldots,n,$ are the corresponding face maps.
\vspace{0.5cm}\\Now let $S\in C^{\infty}(G^{(2)})$ be a cocycle, ie. $\delta^*S=0\,,$ and let $d\mu$ be a Haar measure on $G\,.$ For $f,g\in C_c^{\infty}(G^{(1)})\,,$ their twisted convolution is denoted $f\ast g\in C_c^{\infty}(G^{(1)})\,,$ and is defined by
\vspace{0.05cm}\\
\begin{equation}\label{conv}
   (f\ast g)(\gamma)=\int_{t(\gamma')=t(\gamma)} f(\gamma')g(\gamma'^{-1}\gamma)e^{S[\gamma',\gamma'^{-1}\gamma]}\,d\mu(\gamma')\,.
\end{equation}
Let's observe that 
\begin{equation}
    \{\gamma'\in G^{(1)}:t(\gamma')=t(\gamma)\}\cong \{(\gamma',\gamma'')\in G^{(2)}:\gamma'\gamma''=\gamma\}\,.
\end{equation}
\vspace{0.05cm}\\
Therefore, using the face maps $d_0(\gamma',\gamma'')=\gamma''\,,d_1(\gamma',\gamma'')=\gamma'\gamma''\,,d_2(\gamma',\gamma'')=\gamma'\,,$ we can rewrite \ref{conv} in the following way: 
\vspace{0.05cm}\\
\begin{align}
   (f\ast g)(\gamma)=\int_{\begin{subarray}{l}(\gamma',\gamma'')\in G^{(2)}\\d_1(\gamma',\gamma'')=g\end{subarray}} f(d_2(\gamma',\gamma''))g(d_0(\gamma',\gamma''))e^{S[\gamma',\gamma'']}\,d\mu(\gamma',\gamma'')\,.
\end{align}
\vspace{0.05cm}\\
The advantage of rewriting the twisted convolution this way is that the above formula makes sense on higher groupoids. 
\vspace{0.5cm}\\Let us \textit{emphasize} that, under the usual interpretation of the twisted convolution algebra, we are really just taking the convolution of two sections of a nontrivial multiplicative line bundle (ie. a line bundle over the space of arrows with a compatible multiplication). The total space is given by $G^{(1)}\times\mathbb{C}\,,$ and for a composable pair $(g_2,g_1)\,,$ the multiplication is given by
\vspace{0.05cm}\\
\begin{equation}
    (g_2,\lambda)\cdot (g_1,\beta)=(g_2g_1,\,\lambda\beta \,e^{S[g_2,g_1]})\,.
\end{equation}
\vspace{0.05cm}\\
The trivial multiplicative line bundle corresponds to the zero cocycle, and its twisted convolution algebra just gives the usual convolution algebra. Note that, our construction only gives multiplicative line bundles whose underlying line bundle itself is trivial — it is the multiplication that is nontrivial. However, in general the underlying line bundles need not be trivial, but we won't be needing them; see \ref{appen} for more.
\paragraph{Geometric Quantization}$\,$\\
\vspace{0.05cm}\\Now we will recall what a geometric quantization of a Poisson manifold is, which is used to construct a $C^*$-algebra. Essentially, it is obtained by geometrically quantizing the symplectic groupoid. For full details, see \cite{eli} (also see \cite{weinstein}):
\vspace{0.05cm}\\
\begin{definition}
A geometric quantization of a Poisson manifold $(M,\Pi)$ is given by the following data:
\begin{enumerate}
    \item A symplectic groupoid $(G,\omega)$ integrating $T^*M$ (with the Lie algebroid structure induced by $\Pi\,$).
    \item A Hermitian line bundle with connection $(L,\nabla)\to G^{(1)}$ which prequantizes $\omega\,,$ together with a compatible multiplication on $L\,.$
    \item A Lagrangian polarization of $G^{(1)}$ given by the fibers of a fibration $f:G\to H,$ for some Lie groupoid $H\,,$ such that there exists a multiplicative line bundle $L_H\to H$ with $L\cong f^*L_H\,.$\footnote{Here we are assuming there are no Bohr-Sommerfeld conditions, but this can be done more generally, see \cite{eli}.}
    \item The $C^*$-algebra is then the twisted convolution algebra of $L_H\to H\,.$ 
\end{enumerate}
\end{definition}$\,$
\vspace{0.05cm}\\
One shortcoming of this procedure is that it doesn't give any hint as to how to produce a quantization map, or if one even exists. Another issue is that polarizations are often difficult to find, and sometimes don't exist. However, for the geometric quantization using higher groupoids  we don't need a polarization. Furthermore, the geometric quantization is completely canonical and comes with a quantization map.
\subsubsection{Twisted Convolution Algebra of Higher Integrations of $(M,\Pi)$}$\,$
\vspace{0.05cm}\\
Before beginning this section let us remark that, in the context of path integrals, we will use a physics level of rigour. Now, given a Poisson manifold $(M,\Pi)\,,$ the main simplicial spaces we will be concerned with in this paper are denoted 
\vspace{0.05cm}\\
\begin{equation}
\Pi_1(T^*M)\,,\Pi_2(T^*M)\,,\Pi_{\infty}(T^*M)\,.
\end{equation}
\vspace{0.05cm}\\
The first two are a Lie 1-groupoid (the source simply connected integration of $\Pi$) and a Lie 2-groupoid, respectively, while the third is a simplicial set, however it is not known if it is an $\infty$-groupoid (though in some cases it is, see \cite{getzler}, \cite{andre}). Nevertheless, they are related by taking quotients, ie. 
$\Pi_1(T^*M)$ is a quotient of $\Pi_2(T^*M)\,,$ which is a quotient of $\Pi_\infty(T^*M)$ (see \cite{zhuc}). 
\vspace{0.05cm}\\
\begin{definition}(see \cite{zhuc})
$\Pi_{\infty}(T^*M)$ is a simplicial set defined by
\vspace{0.05cm}\\
\begin{equation}
  \Pi_{\infty}^{(n)}(T^*M)=\{f\in \text{Hom}(\,T|\Delta^n|,\, T^*M\,): f\vert_{\text{vertices}}=0\}\,,
\end{equation}
\vspace{0.05cm}\\
where $|\Delta^n|$ is the standard $n$-simplex (ie. the space — the simplicial set will be denoted by $\Delta^n$). We will denote points in $\Pi_{\infty}^{(1)}(T^*M)$ by $\gamma\,,$ and points in $\Pi_{\infty}^{(2)}(T^*M)$
by $X\,.$
\end{definition}
There is a canonical 2-cocycle on $\Pi_{\infty}(T^*M)\,,$ which for $X\in \Pi^{(2)}_{\infty}(T^*M)$ is given by
\begin{equation}
S[X]:=\frac{i}{\hbar}\int_{|\Delta^2|}X^*\Pi\,.
\end{equation}Now consider two maps 
\begin{equation*}
    s_1\,,s_2: \Pi_{\infty}^{(1)}(T^*M)\to \mathbb{C}\,.
\end{equation*}
Their twisted convolution, denoted 
\begin{equation*}
s_1\ast s_2:\Pi_{\infty}^{(1)}(T^*M)\to \mathbb{C}\,,
\end{equation*}
is given by
\begin{equation}
    (s_1\ast s_2)(\gamma)=\int_{\begin{subarray}{l}X\in \Pi^{(2)}_{\infty}(T^*M)\\d_1(X)=\gamma\end{subarray}}s_1(d_2(X))s_2(d_0(X))e^{\frac{i}{\hbar}S[X]}\,d\mu(X)\,.
\end{equation}\vspace{0.05cm}\\
Here, $d\mu$ is an invariant (or Haar) measure, but we won't go into details here. \textit{Again}, really what we are doing is thinking of $s_1,s_2$ as taking values in the multiplicative line bundle associated to the cocycle $S,$ which we will discuss in more detail later.
\subsubsection{The Product}
Now given $f\in C^{\infty}(M)\,,$ we can lift it to a map $\hat{f}:\Pi^{(1)}_{\infty}(T^*M)\to \mathbb{C}\,,$ 
defined by
\vspace{0.05cm}\\
\begin{equation}
    \hat{f}(\gamma)=(f\circ \pi)(\gamma(1/2))\,,
\end{equation}
\vspace{0.05cm}\\
where $\pi:T^*M\to M$ is the projection. We think of $\hat{f}$ as being the observable associated to $f\,.$
We then get a product on $C^{\infty}(M)$ induced by the twisted convolution algebra, which for $m\in M$ is given by
\vspace{0.05cm}\\
\begin{equation}\label{equ}
    (f\star g)(m)=\int_{\begin{subarray}{l}X\in \Pi^{(2)}_{\infty}(T^*M)\\d_1(X)=m\end{subarray}}\hat{f}(d_2(X))\hat{g}(d_0(X))e^{\frac{i}{\hbar}S[X]}\,d\mu(X)\,.
\end{equation}
\vspace{0.05cm}\\
On the right side of the above equation, we think of $m$ as being the constant algebroid path mapping to $m\,.$ This product agrees with \ref{mystarr}.
\vspace{0.5cm}\\Right now \ref{equ} does not make sense as an integral over $\Pi_1(T^*M)$ do to a lack of gauge invariance. In particular, the observables $\hat{f}\,,\hat{g}$ depend on $\gamma(1/2)\,,$ so they are in no way gauge invariant (in the way that we want). However, we can rewrite \ref{equ} in such a way that the observables only depend on $\gamma(0),\,\gamma(1)\,,$ giving us an integral that makes sense over $\Pi_1(T^*M)\,\footnote{There are some subtleties involving gauge invariance which we will discuss later.}.$ The obvious guess is that the integral over $\Pi_1(T^*M)$ will involve its convolution, but surprisingly it does not.
\vspace{0.5cm}\\Let us remark here that everything we've done makes sense using $\Pi_2(T^*M)$ and all of the formulas are completely unchanged. The only difference is that the phase space gets reduced, ie. $\Pi^{(2)}_2(T^*M)$ is a quotient of $\Pi^{(2)}_{\infty}(T^*M)\,,$ where $X,\,X'\in \Pi^{(2)}_{\infty}(T^*M)$ map to the same point in $\Pi^{(2)}_2(T^*M)$ if they agree on the boundary and are Lie algebroid homotopy equivalent, relative to the boundary. For this reason we are going to use $\Pi_2$ and $\Pi_{\infty}$ interchangeably.
\subsubsection{Groupoid-Valued Sigma Models}$\,$
\vspace{0.05cm}\\
We will briefly discuss the Poisson sigma model here and the higher Poisson sigma models which determine the star products $\star\,.$
\paragraph{Poisson Sigma Model}$\,$
\vspace{0.5cm}\\The Poisson sigma models studied in \cite{catt} has the following ingredients:
\begin{enumerate}
    \item The parameter (or worldsheet) and target spaces are $D,M,$ respectively, where $D$ is a disk with three marked points and $(M,\Pi)$ is a Poisson manifold.
    \item The fields are vector bundle morphisms $(X,\eta):TD\to T^*M$ covering $X:D\to M\,,$ such that the pullback of $\eta$ to $\partial D$ takes values in the zero section.
    \item The action of a field is given by 
    $S[X,\eta]=\int_D \eta\wedge dX+\frac{1}{2}\Pi\vert_X(\eta\,,\eta)\,.$
    \item For $f,g\in C^{\infty}(M)\,,$ we define $
 (f\star g)(m)=\int_{X(\infty)=m} f(X(1))g(X(0))\,e^{\frac{i}{h}S[X,\eta]}\,dX\,d\eta\,.$
\end{enumerate}$\,$
\vspace{0.05cm}\\
\paragraph{``Sigma Model" of the Lie Algebroid}$\,$
\vspace{0.5cm}\\Closely related to the Poisson sigma model is another ``sigma model" which emphasizes the Lie algebroid structure of $\pi:T^*M\to M\,.$ Its ingredients are:
\begin{enumerate}
    \item The parameter and target spaces are $|\Delta^2|, T^*M,$ respectively.
    \item The fields are Lie algebroid morphisms $X:T|\Delta^2|\to T^*M\,.$
     \item The action of a field is given by 
    $S[X]=\int_{|\Delta^2|} X^*\Pi$
    \item For $f,g\in C^{\infty}(M)\,,$ we use the twisted convolution algebra to define 
    \begin{equation}\label{inttt}
        (f\star g)(m)=\int_{\pi(X(\infty))=m}f(\pi(X(1)))\,g(\pi(X(0)))\,e^{\frac{i}{\hbar}S[X]}\,dX\,,
        \end{equation}
        where we identify $|\Delta^2|$ with the disk $D$ with three marked points.
\end{enumerate}$\,$
\vspace{0.05cm}\\
\paragraph{Sigma Model of $\mathbf{\Pi_2(T^*M)}$}$\,$
\vspace{0.5cm}\\Looking at the ``sigma model" of the Lie algebroid, we see that the fields are points in $\Pi^{(2)}_2(T^*M)\,,$ and by the Yoneda lemma, \\
\begin{equation}
\Pi^{(2)}_2(T^*M)\cong \text{Hom}(\Delta^2,\Pi_2(T^*M))\,.
\end{equation}$\,$\\
Here, $\Delta^n$ is the standard $n$-simplex (ie. the simplial set). We think of this as a nerve of the following category:
\begin{itemize}
    \item The objects are points in $\{0,1,\ldots,n\}\,,$
    \item For each pair of objects $j,i\,,$ if $j\ge i$ then there is exactly one morphism $j\leftarrow i\,,$ and if $j<i$ there are no morphisms $j\leftarrow i\,.$ For $j\ge i\,,$ we denote the unique morphism by $(j,i)\,.$
\end{itemize}$\,$\\
We will now discuss a sigma model (or something analogous) og $\Pi_2(T^*M)\,.$ Its ingredients are:\\
\begin{enumerate}
    \item The parameter and target spaces are $\Delta^2,\, \Pi_2(T^*M),$ respectively.
    \item  The fields are morphisms of simplicial spaces $F:\Delta^2\to \Pi_2(T^*M)\,,$ which in degree $n$ will be denoted by $F^{(n)}\,.$
    \item The Lagrangian density is, for $X\in\Pi_2^{(2)}(T^*M)\,,$ $S[X]=\int_{|\Delta^2|} X^*\Pi\,.$
    \item The action of a field $F$ is given by 
    \begin{equation}
        \sum_{(b,a)\in\Delta^2\text{ are a  composable pair}}S[F^{(2)}(b,a)]=S[F^{(2)}((2,1),(1,0))]\,.
        \end{equation}
        \item For $f,g\in C^{\infty}(M)\,,$ we define 
        \begin{equation}\label{inttttt}
            (f\star g)(m)=\int_{\begin{subarray}{l}F\in\text{Hom}(\Delta^2,\Pi_2(T^*M)),\\F(1,1)=m\end{subarray}} dF_m\,e^{\frac{i}{\hbar}S(F(2,1),F(0,1))}\,f(F(2,2))g(F(0,0)) \,.
              \end{equation}
              Here $dF_m$ is a measure on the space $\{F\in\text{Hom}(\Delta^2,\Pi_2(T^*M)): F(1,1)=m\}\,.$ Again, we think of $m$ as being the constant algebroid path mapping to $m\,.$
\end{enumerate}
A crucial point is that the products defined by \ref{inttt},\ref{inttttt} are the same. Later in this paper we will explain why.
\vspace{0.05cm}\\
\paragraph{Sigma Model of $\mathbf{\Pi_1(T^*M)}$}$\,$
\vspace{0.5cm}\\While it's not clear how to generalize the first two sigma model constructions to groupoids, the third one does generalize, assuming certain descent conditions. For now, we will assume that the cocycle 
\vspace{0.05cm}\\
\begin{equation*}
    S[X]=\int_{|\Delta^2|}X^*\Pi
    \end{equation*}
    \vspace{0.05cm}\\
descends to a cocycle on $\Pi_1(T^*M);$ we will explain precisely what this means later, but it is enough that, for all Lie algebroid morphisms $X:TS^2\to T^*M\,,$
\vspace{0.05cm}\\\begin{equation}
    \int_{S^2}X^*\Pi=0\,.
    \end{equation}
    \vspace{0.05cm}\\
This happens, for example, whenever the Poisson structure is exact in Poisson cohomology, or whenever the second homotopy group of each source fiber (or symplectic leaf) is finite, eg. Lie-Poisson structures and symplectic tori, respectively. Though this isn't the case for the symplectic 2-sphere, for example, we will see that our construction can still be applied to it.
\vspace{0.5cm}\\We will now discuss the groupoid-valued sigma model (groupoid-valued in the sense that the target space is a groupoid). Its ingredients are:\\
\begin{enumerate}
    \item The parameter and target spaces are $\Delta^2,\, \Pi_1(T^*M),$ respectively.
    \item  The fields are morphisms of simplicial spaces $F:\Delta^2\to \Pi_1(T^*M)\,,$ which in degree $n$ will be denoted by $F^{(n)}\,.$
    \item The Lagrangian density is given by a real-valued 2-cocycle $S$ on $\Pi_1(T^*M)\,,$ such that:
    \begin{itemize}
        \item $VE(S)=\Pi$ (ie. the van Est map applied to $S$ equals the 2-cocycle associated with the Poisson structure),
        \item  $S$ is antisymmetric, in the sense that:
    \begin{equation}
        S(g_2,g_1)=-S(g_1^{-1},g_2^{-1})\,.
    \end{equation}
    \end{itemize} 
    \item The action of a field is given by 
    \begin{equation}
        \sum_{(b,a)\in\Delta^2\text{ are a  composable pair}}S[F^{(2)}(b,a)]=S[F^{(2)}((2,1),(1,0))]\,.
        \end{equation}
        \item For $f,g\in C_c^{\infty}(M)\,,$ we define a (possibly nonassociative) product by
        \begin{equation}\label{intttt}
            (f\star g)(m)=\int_{\begin{subarray}{l}F\in\text{Hom}(\Delta^2,G),\\F(1,1)=m\end{subarray}} dF_m\,e^{\frac{i}{\hbar}S(F(2,1),F(0,1))}\,f(F(2,2))g(F(0,0)) \,.
              \end{equation}
              Here $dF_m$ is a measure on the space $\{F\in\text{Hom}(\Delta^2,G): F(1,1)=m\}\,.$ On the right, we think of $m$ as being the identity morphism at $m\,.$
              \item Addtionally, one may want to equip the Poisson manifold with a volume form in order to compute all correlation functions. Of course, the symplectic leaves already come with a (Hamiltonian invariant) volume form.
\end{enumerate}$\,$\\
Later in this paper we will explain how to find the measure $dF_m\,,$ but it is an actual measure. Note that, if one wanted to avoid simplicial sets that aren't groupoids, one could replace $\Delta^2$ with $\text{Pair}\{0,1,2\}\rightrightarrows\{0,1,2\}$ and essentially everything remains the same. 
\vspace{0.5cm}\\Graphically, using the identification of composable pairs with morphisms $\Delta^2\to \Pi_1(T^*M)\,,$ \ref{intttt} can be represented in the following way:
\[
\begin{tikzpicture}
\tikzstyle arrowstyle=[scale=2]  \tikzstyle directed=[postaction={decorate,decoration={markings,mark=at position .5 with{\arrow[arrowstyle]{stealth}}}}]
\foreach \i in {0, 180, 270}{
\filldraw (\i:1.25)circle(0.05);
}
\draw[dashed, directed] (0:1.25)arc[start angle=0, end angle=180, x radius=1.25cm, y radius=1.25cm] node[midway, above] {$\gamma_2\gamma_1\;\;\,$};
\draw[directed] (270:1.25)arc[start angle=-90, end angle=-180, x radius=1.25cm, y radius=1.25cm] node[midway, left] {$\gamma_2$};
\draw[directed] (0:1.25)arc[start angle=0, end angle=-90, x radius=1.25cm, y radius=1.25cm] node[midway, right] {$\;\,\gamma_1$};

      \path
      (-2,0) node[] {$\displaystyle\bigint$}
      (-3.5, 0) node[] {$(f\star g)(m)\;=$}
      (-1.25, 0) node[left] {$f$}
      (1.25, 0) node[right] {$g$}
      (0, -1.25) node[below] {$m$}
      (0, 0) node[] {$e^{\frac{i}{\hbar}S[\gamma_2,\gamma_1]}$}
      (-2.25, -1.2) node[] {$\gamma_1,\gamma_2$}
      
      ;
\end{tikzpicture}
\]
We see that, in the spirit of Feynman, $f\star g$ can be interpreted of as a sum over all paths, where the paths are given by arrows in the Lie groupoid. Using the formulation of the star product given in \cite{bon}, the paths are arrows in the corresponding Lie 2-groupoid. Now even though this sigma model is finite dimensional (in the sense that the space of fields is finite dimensional), there are still sometimes infinities to deal with due to gauge invariance. However, these infinities can often be overcome non-perturbatively, as we will see.\\
\vspace{0.05cm}\\One way of interpreting $(f\star g)(m)$ is as the groupoid version of the product: any two objects connected by an arrow are isomorphic, so $(f\star g)(m)$ ``averages" $fg$ over all objects isomorphic to $m\,.$\\
\vspace{0.05cm}\\
From the above definitions it is clear that we are defining ``products" via integrals with kernels involving complex exponentials, and this is quite typical among formulas for strict deformation quantizations/convergent products. See \cite{rieffel} for several examples, as well as \cite{Pierre}; for a survey of convergent products, see \cite{Waldmann}.
\vspace{0.05cm}\\
\begin{remark}\label{asso}Let us \textit{emphasize} here that we are not claiming that these products are associative, but we will see several examples in which they are. However, even when the product isn't associative, we still obtain a $C^*$-algebra that one would hope forms a strict quantization (as discussed in \ref{explain}, after dealing with any infinities).
\end{remark}$\,$
\section{Essential Example}\label{fund}$\,$
\vspace{0.05cm}\\
Now let's flesh out an example of a groupoid-valued sigma model. We will begin with the most basic nontrivial example; we will start from the Poisson sigma model of symplectic $\mathbb{R}^2$ and derive various formulations of quantum mechanics.
\\\\The point is to do this example in such a way that the constructions generalize to other Poisson manifolds. Points 1-7 involve constructing the Moyal product starting from the sigma model (see \cite{rieffel} for an alternative derivation); points 8-14 involve constructing irreducible representations of the $C^*$-algebra. This involves a bit of $C^*$-algebra theory (which we will explain), for references see \cite{Gleason}, \cite{Murphy};  points 15-16 involve a phase space formulation of quantum mechanics, dating back to Moyal, Weyl, Wigner, see \cite{Curtright}, \cite{Moyal}, \cite{Wigner}. We will cover a bit more detail in example \ref{symp}, but for now:
\vspace{0.05cm}\\
\begin{exmp}
Consider the constant symplectic structure $(\mathbb{R}^2,\omega^{-1}=\partial_p\wedge\partial_q)\,.$ 
\begin{enumerate}
    \item The parameter and target spaces are  $\Delta^2,\;\Pi_1(T\mathbb{R}^2)\cong \mathbb{R}^2\ltimes\mathbb{R}^2\rightrightarrows \mathbb{R}^2\,,$ respectively (we are using the identification $T\mathbb{R}^2\cong T^*\mathbb{R}^2$ via the symplectic form). That is, the target space is the action groupoid of $\mathbb{R}^2$ acting on itself by addition. We may denote this groupoid simply by $\mathbb{R}^2\ltimes\mathbb{R}^2\,.$
    \item The fields are morphisms of simplicial spaces $F:\Delta^2\to \mathbb{R}^2\ltimes\mathbb{R}^2\,,$ which in degree $n$ will be denoted by $F^{(n)}\,.$ These fields can be identified with composable pairs in $\mathbb{R}^2\ltimes\mathbb{R}^2\,,$ via the map 
    \begin{equation*}
        F\mapsto F^{(2)}((2,1),(1,0))\,.
        \end{equation*}
        \item The cocycle is given by 
        \begin{equation}
            S[(a',b',a+p,b+q),(a,b,p,q)]=\frac{a'b-b'a}{2}\,.
        \end{equation}
             This is the area of a triangle in the plane with vertices at $(p,q),\,(a+p,b+q),\,(a'+a+p,b'+b+q)\,.$ 
        \vspace{0.05cm}\\This cocycle can be obtained by gauge fixing the Poisson sigma model, in the following way: first, the fields in this Poisson sigma model are equivalently maps $X:|\Delta^2|\to \mathbb{R}^2\,.$ Now given any composable pair in $\mathbb{R}^2\ltimes\mathbb{R}^2$ we can canonically form the boundary of a 2-simplex $\partial X$ by joining the vertices with straight lines. No matter how we fill in this 2-simplex we have that 
        \begin{equation*}
            \int_{|\Delta^2|}X^*\omega=\text{symplectic area of the corresponding triangle in $\mathbb{R}^2$}\,,
        \end{equation*}
        which is exactly what our cocycle is.
        \item Using the identification of the fields with composable pairs, we get the following identification:
        \begin{equation*}
            \{F\in\text{Hom}(\Delta^2,\mathbb{R}^2\ltimes\mathbb{R}^2): F(1,1)=(p,q)\}\cong\{\big((a',b',p,q),(a,b,p-a,q-b)\big)\}\subset(\mathbb{R}^2\ltimes\mathbb{R}^2)^{(2)}\,.
        \end{equation*}
        Therefore, the former space is identified with $\mathbb{R}^4\,,$ with coordinates $(a',b',a,b)\,.$
        \item There is a natural Haar measure on the source and target fibers, given by $t^*(dp\wedge dq)=da'\wedge db'\,,\,s^*(dp\wedge dq)=da\wedge db\,,$ respectively. Using the previous identification, the measure is given by 
        \begin{equation*}
            dF_{(p,q)}=da'\wedge db'\wedge da\wedge db\,.
        \end{equation*}
        \item For $f,g\in C_c^{\infty}(\mathbb{R}^2)\,,$ we have that
        \begin{equation*}
            (f\star g)(p,q)=\frac{1}{(4\pi\hbar)^2}\int_{\mathbb{R}^4}f(a'+p,b'+q)g(p-a,q-b)e^{\frac{i}{2\hbar}(a'b-b'a)}\,da'\,db'\,da\,db\,.
        \end{equation*}
        Here, $1/(4\pi\hbar)^2$ is a normalization constant, determined by $1\star 1=1$ ($1\star 1$ should be considered as an improper integral). Therefore, the measure we are really interested in is
        \begin{equation*}
            \frac{1}{(4\pi\hbar)^2}da'\wedge db'\wedge da\wedge db\,.
        \end{equation*} \\
        Changing variables, we can rewrite this as:
        \begin{align}\label{moy}
          &(f\star g)(p,q)=\nonumber\\&\boxed{\frac{1}{(4\pi\hbar)^2}\int_{\mathbb{R}^4}f(p'',q'')g(p',q')e^{\frac{i}{2\hbar}[(p''-p)(q-q')-(q''-q)(p-p')]}\,dp''\,dq''\,dp'\,dq'}\,.  
        \end{align}
        This defines an \textit{associative product} (see \cite{Zachos}).
    \item This is the integral form of the Moyal product (see \cite{baker}, where it was originally introduced). Its asymptotic expansion in $\hbar\,,$ about $\hbar=0\,,$ gives the usual perturbative version of the Moyal product (see \cite{kontsevich}):
    \begin{equation*}
        f\star g\sim m\circ e^{i\frac{\hbar}{2}\partial_p\wedge\partial_q}(f\otimes g)\;\;\;(\hbar\to 0)\,,
    \end{equation*}
    here $m$ denotes multiplication. This can be obtained by using the stationary phase approximation, differentiating in $\hbar$ and integrating by parts (see example \ref{good}).\\\\
    \item There is a natural $C^*$-algebra and representation associated with $\star\,,$ as discussed in \ref{explain}. It is given by the completion, with respect to the operator norm, of the left $\star$-action of $C^{\infty}_c(\mathbb{R}^2)$ on $L^2(\mathbb{R}^2,dp\wedge dq)\,.$ We have that $\overline{f\star g}=\bar{g}\star\bar{f}\,,$ and indeed complex conjugation gives the involution of $C_c^{\infty}(\mathbb{R}^2)$ inside this $C^*$-algebra.
    \item Now in order to do quantum mechanics, one just needs to pick a state, which are given by normalized functions in our Hilbert space. However, in quantum mechanics one is usually interested in states that live in irreducible representations, and this one isn't. However, $C^*$-algebras come with a notion of pure states, which are the states which live in irreducible representations.\footnote{See \cite{Gleason}, \cite{Murphy}.}
    \begin{definition}
        A state on a unital $C^*$-algebra $\mathcal{A}$ is a linear functional $\rho:\mathcal{A}\to\mathbb{C}$ such that
        \begin{itemize}
            \item $\rho(a^*a)\ge 0$ for all $a\in\mathcal{A}\,,$
            \item $\rho(1)=1\,.$
        \end{itemize}
   States are closed under convex combinations. A state is called pure if it isn't a convex combination of distinct states, otherwise the state is called mixed.
    \end{definition}
    Now associated to a representation of a $C^*$-algebra on a Hilbert space are vector states, which is how states are most commonly thought of in quantum mechanics:
    \begin{definition}
        Let $\pi:\mathcal{A}\to \mathcal{B}(\mathcal{H})$ be a $C^*$-algebra representation on $\mathcal{H}\,.$ A vector state is a state of the form $a\mapsto \langle \psi,a\psi\rangle\,,$ for some $\psi\in \mathcal{H}\,.$
    \end{definition}
   By the GNS construction, every state can be represented by a vector state, in some representation. Now, associated to every $\psi\in\mathcal{H}$ is a representation of $\mathcal{A}$ (the intersection of all subrepresentations containing $\psi$). Now to relate pure states to irreducibility, we have the following theorem:
    \begin{theorem}\label{red}(see \cite{Gleason}, \cite{Murphy})
        The vector state associated to $\psi$ is pure if and only if the associated representation is irreducible. 
    \end{theorem}
    The above theorem implies that all states in an irreducible representation are pure.\\
    \vspace{0.05cm}\\
    Now in our representation there are many vector states, and some of them define the same state. Therefore, our state space really consists of a quotient of the normalized functions in $L^2(\mathbb{R}^2)\,,$ where $\psi\sim \psi'$ if they define the same state. Furthermore, this space contains a subspace of pure states. \Cref{red} tells us that in order to find pure states in our Hilbert space it is enough to look for irreducible representations (in fact, by the Stone-von Neumann theorem, since our Hilbert space does contain an irreducible representation, it contains all pure states). 
    \item One such irreducible representation is given as follows: the left action given by \ref{moy} preserves functions of the form 
    \begin{equation}\label{form}
        \psi(p+iq)e^{-\frac{(p^2+q^2)}{2\hbar}}\,,
    \end{equation}
    where $\psi(p+iq)$ is holomorphic in $z=p+iq\,.$ Therefore, we should consider the closed subspace of our Hilbert space generated by functions of the form \ref{form}. The functions $z,\,\bar{z}$ act as follows:
    \begin{align}
        &z\star  \psi(z)e^{-\frac{|z|^2}{2\hbar}}=z\psi(z)e^{-\frac{|z|^2}{2\hbar}}\,,
        &\bar{z}\star  \psi(z)e^{-\frac{|z|^2}{2\hbar}}=\Big[\hbar\frac{\partial}{\partial z}\psi(z)\Big]e^{-{\frac{|z|^2}{2\hbar}}}\,.
    \end{align}
    This is Segal–Bargmann space, one of the common representations used in quantum mechanics, see \cite{hall}.
        \item We can generate submodules (and subrepresentations) for the algebra in the following way: let $a,b,c,d\in\mathbb{C}$ be such that $ad-bc=1\,.$ Consider $A=aq+bp,\,C=cq+dp$ (then \{A,C\}=1).\footnote{These are canonical coordinates, in the classical mechanics sense.} Functions of the form
        \vspace{0.05cm}\\
    \begin{equation}
        \psi(p,q)=e^{\frac{i}{2\hbar}(aq+bp)(cq+dp)}\psi(cq+qp)
    \end{equation}
    \vspace{0.05cm}\\
form a submodule for the algebra. Here $\psi(cq+dp)$ should be interpreted as functions which satisfy\footnote{This is the one place where, in a sense, we are using a polarization.}
\vspace{0.05cm}\\
\begin{equation}
    \frac{1}{2}\Big(d\frac{\partial}{\partial q}-c\frac{\partial}{\partial p}\Big)\psi=0\,,
\end{equation}

and $aq+bp, cq+dp$ act by
\begin{align}
    &(aq+bp)\star \Big[e^{\frac{i}{2\hbar}(aq+bp)(cq+dp)}\psi(cq+qp)\Big]=e^{\frac{i}{2\hbar}(aq+bp)(cq+dp)}\Big[\frac{2\hbar}{i}\psi'(cq+dp)\Big]\,,
    \\&(cq+dp)\star \Big[e^{\frac{i}{2\hbar}(aq+bp)(cq+dp)}\psi(cq+qp)\Big]=e^{\frac{i}{2\hbar}(aq+bp)(cq+dp)}\big[2(cq+bp)\psi(cq+qp)\big]\,.
\end{align}
\vspace{0.05cm}\\
Of course, these submodules come in pairs since functions of the form 
\vspace{0.05cm}\\
\begin{equation}
        \psi(p,q)=e^{-\frac{i}{2\hbar}(aq+bp)(cq+dp)}\psi(aq+bp)
    \end{equation}
    \vspace{0.05cm}\\
    must also form a submodule (since if $ad-bc=1,$ then $c(-b)-d(-a)=1$). This family of submodules therefore (essentially) contains the usual position and momentum representations, as well as the Segal-Bargmann represenation. In the case that $a,b\in\mathbb{R},\,$ $|\psi(aq+bp)|^2$ is the probability density of obtaining $aq+bp$ when performing the corresponding measurement.
        \item Now $(f\star g)(p,q)$ is supposed to be a 3-point function (see \cite{catt}). In general, the expectation value of a function of the fields $\Phi$ is given by
        \begin{align}\label{exp}
          \langle\Phi\rangle=\frac{1}{(4\pi\hbar)^2}\int_{\mathbb{R}^6}\Phi(p'',q'',p',q',\mathfrak{p},\mathfrak{q})e^{\frac{i}{2\hbar}[(p''-p)(q-q')-(q''-q)(p-p')]}\,dp''\,dq''\,dp'\,dq'\,d\mathfrak{p}\,d\mathfrak{q}\,.  
        \end{align}
        Then what we have is 
        \begin{align}
           &(f\star g)(p,q)=\nonumber
           \\& \frac{1}{(4\pi\hbar)^2}\int_{\mathbb{R}^6}f(p'',q'')g(p',q')\delta_{(p,q)}(\mathfrak{p},\mathfrak{q})e^{\frac{i}{2\hbar}[(p''-p)(q-q')-(q''-q)(p-p')]}\,dp''\,dq''\,dp'\,dq'\,d\mathfrak{p}\,d\mathfrak{q}\,.  
        \end{align}
        \vspace{0.05cm}\\
        Here $\delta_{(p,q)}(\mathfrak{p},\mathfrak{q})$ corresponds to the observable $F\mapsto \delta_{(p,q)}(F(1,1))\,.$ Therefore, we indeed have the following equality:
        \vspace{0.05cm}\\
        \begin{equation}
           (f\star g)(p,q)= \langle\, f(F(2,2))g(F(0,0))\delta_{(p,q)}(F(1,1))\,\rangle\,.
            \end{equation}
            \item Given a function $f(p,q),$ we can compute the expectation value in two different ways: one is obtained by computing the expectation value of $\hat{f}$ via \ref{exp} (we can define $\hat{f}$ by $\hat{f}(F)=f(F(j,j))$ for any fixed $j=0,1,2,$ they all have the same expectation) and the other is by computing the expectation value as an operator in our quantum system. We have the following relation:
      \begin{equation}
          \langle \hat{f}\rangle =\langle 1,f\star 1\rangle=\int_{\mathbb{R}^2}f(p,q)\,dp\,dq\,.
      \end{equation}
      Therefore the computations agree — note that, we are computing the expectation values in the unnormalizable state $h(p,q)=1,$ but we can always use a normaizable one. In fact, the state we are using is the maximal mixed state,\footnote{A maximally mixed state is a state whose density matrix is proportional to the identity. One could show that formally, the density matrix is given by $\rho(p,q)=\frac{(2\pi\hbar)^2}{\int_{\mathbb{R}^2}dp\,dq}\,,$ ie. $\langle f\rangle=\text{tr}(\rho\star f)$, for any $f\,.$} which is essentially a state in which a hypothetical physicist has no information about the wave function.
      \item Let us remark one interesting (well-known) property: we have
\begin{equation}
\int_{\mathbb{R}^2}f\star g\,dp\,dq=\int_{\mathbb{R}^2}fg\,dp\,dq\,.
\end{equation}
This is just a consequence of the fact that $L^2(\mathbb{R}^2)$ forms a $C^*$-algebra representation, ie.
\begin{align}
   & \langle f\star g\rangle=\langle 1,  f\star g\star 1\rangle= \langle \overline{f}\star 1,g\star 1\rangle
   =\langle\overline{f},g\rangle\,.
\end{align}
\item Essentially, we've constructed a phase space operator formulation of quantum mechanics. Of course, there is a phase space quasiprobability distributon formulation of quantum mechanics that dates back to Moyal, Weyl, Wigner, and others.\footnote{See \cite{Curtright}, \cite{Moyal}, \cite{Wigner}.} These formulations are related through the Riesz representation theorem, as follows: let $\rho$ be a state on our $C^*$-algebra. By the Riesz representation theorem, there exists a $P\in L^2(\mathbb{R}^2)$ such that, for any observable $h(p,q)\,,$
\begin{equation}
    \rho(h)=\langle W,h\rangle\,.
\end{equation}
It follows that $P$ is real-valued and
\begin{equation}
    \int_{-\infty}^{\infty}W(p,q)\,dp\,dq=1\,,
    \end{equation}
and we say that $W$ is a quasiprobability distribution. The Wigner distribution is then obtained by applying this construction to vector states in usual quantum mechanics. That is, given a wave function $\psi(q)\,,$ the Wigner distribution is a real-valued function $W(p,q)$ integrating to $1,$ such that for any observable $h(p,q)\,,$
\begin{equation}
    \langle h,W\rangle_{L^2(\mathbb{R}^2)} =\langle \psi, h\star\psi\rangle_{L^2(\mathbb{R})}\,.
\end{equation}
Given a Hamiltonian $H(p,q)\,,$ the time evolution of an observable $f$ is as in the Heisenberg picture,\footnote{Though it can also be done using the Schr\"{o}dinger picture.} ie. it is given by Hamilton's equations: 
\begin{equation}
    \frac{df}{dt}=[f,H]:=f\star H-H\star f\,.
\end{equation}
The expectation value of $f_t(p,q)$ is just $\langle f_t,W\rangle\,. $
\end{enumerate} 
\end{exmp}$\,$\\
Now that we've given an overview of this paper we will get into the details and do more examples.
\section{Higher Integrations of $(M,\Pi)$}
In this section we are going to discuss $\Pi_2(T^*M)\,,$ which is a Lie 2-groupoid integrating the Poisson structure $\Pi$ (see \cite{zhuc}). The first step is defining the simplicial set $\Pi_{\infty}(T^*M)\,.$ Ideally, $\Pi_{\infty}(T^*M)$ would turn out to be an $\infty$-groupoid, but in general this isn't known (though in some cases it is, see \cite{getzler}, \cite{andre}).\\
\vspace{0.05cm}\\
Let us remind the reader that in the context of path integrals we will use a physics level of rigour.\\
\vspace{0.05cm}\\Recall that $|\Delta^n|$ denotes the standard $n$-simplex (ie. the space), we will use $\Delta^n$ to denote the simplicial set.\\
\vspace{0.05cm}\\
\begin{definition}(see \cite{zhuc})
$\Pi_{\infty}(T^*M)$ is a simplicial set defined by
\vspace{0.05cm}\\
\begin{equation}
  \Pi_{\infty}^{(n)}(T^*M)=\{f\in \text{Hom}(\,T|\Delta^n|,\, T^*M\,): f\vert_{\text{vertices}}=0\}\,.
\end{equation}
\vspace{0.05cm}\\
Note that these are morphisms of Lie algebroids. Therefore,
\vspace{0.05cm}\\
\begin{enumerate}
    \item The objects of $\Pi_{\infty}(T^*M)$ are the points in $M\,,$ denoted $m\,.$
    \item The arrows in $\Pi_{\infty}(T^*M)$ are $A$-paths, denoted $\gamma$ (we may use this notation to denote arrows in groupoids, in general).
    \item The composable pairs\footnote{Technically these aren't just composable pairs of arrows, but we will use this terminology.} in $\Pi_{\infty}(T^*M)$ are given by $\text{Hom}(T|\Delta^2|, T^*M)\,,$ denoted $X\,.$
    \item The composable $n$-tuples in $\Pi_{\infty}(T^*M)$ are given by $\text{Hom}(T|\Delta^n|, T^*M)\,,$ denoted $X^{(n)}$ (most of the time we will use the previous notation for objects, arrows and composable arrows).
\end{enumerate}
\end{definition}$\,$
\vspace{0.05cm}\\
For a composable pair $X\,,$ we think of the three boundary components,  denoted $\gamma_1,\gamma_2,\gamma_3,$ as being arrows. We think of $\gamma_3$ as being the composition of $\gamma_2$ and $\gamma_1$ and we say that $X$ composes to $\gamma_3\,.$ Diagrammatically, we have the following:
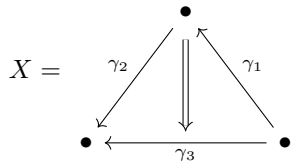
\[X=\begin{tikzcd}
& |[alias=C]| \bullet \arrow[dl,"\gamma_2"'] &  \\[2em]
\bullet   
&  & \bullet  \arrow[ul, "\gamma_1"'] \arrow[ll, "\gamma_3" {name=U}]
\arrow[Leftarrow,from=U,to=C,shorten >=1ex,shorten <=1ex]
\end{tikzcd}\]
\vspace{0.05cm}\\When it comes to this paper $\Pi_{\infty}(T^*M)$ and  $\Pi_{2}(T^*M)$ are mostly interchangeable, but for completion it is nice that it exists as a Lie 2-groupoid (and it gives a partially reduced phase space for the Poisson sigma model). Note that $\Pi_{1}(T^*M)$ refers to the source simply connected symplectic groupoid. \\
\vspace{0.05cm}\\
\begin{definition}(see \cite{zhuc} for details)
$\Pi_2(T^*M)$ is a Lie 2-groupoid, whose space of objects and arrows are the same as $\Pi_{\infty}(T^*M)\,,$ and whose composable pairs are given by $\Pi_{\infty}^{(2)}(T^*M)/\Pi_{\infty}^{(3)}(T^*M)\,,$
\end{definition}$\,$
\vspace{0.05cm}\\
\section{Geometric Quantization}
Here we are going to describe the geometric quantization of $(M,\Pi)$ on $\Pi_{\infty}(T^*M)\,.$ We use the term ``geometric quantization" a bit loosely: rather than having a symplectic form we will have a presymplectic form — it will be a closed 2-form obtained by pulling back a symplectic form to a coisotropic submanifold.\\\\
After constructing the symplectic groupoid, the first step in the geometric quantization of Poisson manifolds normally involves choosing a multiplicative line bundle with a connection whose curvature is the symplectic form, which can equivalently be described as a cocycle in degree 2. On $\Pi_{\infty}(T^*M)$ there is a natural presymplectic form, connection and multiplicative line bundle, however we will mostly only be needing the multiplicative line bundle.\\
\vspace{0.05cm}\\
In several places in this paper we will use homotopies of Lie algebroid morphisms, which is a generalization of homotopies from topology. For references (and more generally for references about constructing the groupoids), see \cite{bojo}, \cite{rui}, \cite{ruif}. We will give the definition:
\begin{definition}(see \cite{bojo})
Let $A_1\to M_1,\,A_2\to M_2$ be Lie algebroids. We say that two morphisms $f,f':A_1\to A_2$ are homotopy equivalent if there exists a morphism $F:A_1\oplus TI\to A_2$ such that the restrictions $F\vert A_1\times \{0\},\,F\vert A_1\times \{1\}$ agree with $f,\,f',$ respectively. Here $I=[0,1]$ and $A_1\oplus TI$ is the direct product Lie algebroid over $A_1\times I\,.$
\end{definition}
\subsection{The Canonical Cocycle}
Before defining the canonical cocycle we will observe that, given a composable $n$-tuple $X^{(n)}:T\Delta^n\to T^*M\,,$ we can pullback $\Pi$ to get a two-form on $\Delta^n$ given by 
\vspace{0.05cm}\\
\begin{equation}
    (v,w)\mapsto \big(X^{(n)*}\Pi\big)(v,w):=\Pi(X^{(n)}(v),X^{(n)}(w))\,.
\end{equation}
\vspace{0.05cm}\\
where $v,w\in T\Delta^n$ are vectors over the same point. Since $X^{(n)}$ is a Lie algebroid morphism, we have that 
\vspace{0.05cm}\\
\begin{equation}\label{closed}
    dX^{(n)*}\Pi=X^{(n)*}d_{\Pi}\Pi=0\,,
\end{equation}
\vspace{0.05cm}\\
where $d_\Pi$ is the Chevalley-Eilenberg differential\footnote{See \ref{appen}.} (note that $d_{\Pi}(\Pi)=[\Pi,\Pi]=0$).
We can then define a map  
\vspace{0.05cm}\\
\begin{equation}\label{myaction}
\Pi_{\infty}^{(2)}(T^*M)\to\mathbb{R}\,,\;X\mapsto S[X]:=\int_{|\Delta^2|}X^*\Pi\,,
\end{equation}
\vspace{0.05cm}\\
where $S$ may be called ``the action".
This map is a 2-cocycle on $\Pi_{\infty}(T^*M)$, and this follows from the fact that, for a composable triple $X^{(3)}:T\Delta^3\to T^*M\,,$
\vspace{0.05cm}\\
\begin{equation}
     (\delta^*S)(X^{(3)})=\int_{\partial\Delta^3}X^{(3)*}\Pi=\int_{\Delta^3}dX^{(3)*}\Pi=0\,,
\end{equation}
\vspace{0.05cm}\\
where we have used Stokes' theorem and \cref{closed}. Now we may exponentiate the action to get a cocycle valued in $S^1$:
\vspace{0.05cm}\\
\begin{definition}
The canonical cocycle on $\Pi_{\infty}(T^*M)$ (or $\Pi_{2}(T^*M)$) is given by $X\mapsto e^{\frac{i}{\hbar}S[X]}\,.$
\end{definition}$\,$
\vspace{0.05cm}\\This cocycle fits into another cocycle with respect to the complex of sheaves\footnote{Technically we should only use this complex of sheaves on the simplicial space in degree higher than 0, and in degree 0, ie. on $M,$ we should put the trivial complex of sheaves $\{e\}\to\{e\}\to\cdots.$ This is discussed in \cite{Lackman} under the term ``truncation".} 
\begin{equation}
    \mathcal{O}^*\xrightarrow[]{\text{dlog}} \Omega^1\xrightarrow[]{d}\Omega^2\xrightarrow[]{d} \cdots,
\end{equation} given by a triple $(e^{\frac{i}{\hbar}S}, \theta,\omega)\,.$\footnote{We won't be needing all of the details.} What this means is that 
\vspace{0.05cm}\\\begin{equation}
d\omega=0,\;\delta^*\omega=0,\;d\theta=\omega,\;\delta^*\theta=i\hbar\,\text{dlog}\,e^{\frac{i}{\hbar}S[X]},\;\delta^*e^{\frac{i}{\hbar}S[X]}=1\,.
\end{equation}
\vspace{0.05cm}\\
Here $\theta,\,\omega$ are a 1-form and 2-form on $\Pi_{\infty}^{(1)}(T^*M),$ which can be described in the following way: for a space $A$ let $PA$ denote the path space, ie. the space of maps $[0,1]\to A.$ Now $\Pi_{\infty}^{(1)}(T^*M)$ is a (infinite dimensional) coisotropic submanifold of $P(T^*M),$ where the latter comes with a (Liouville) 1-form and (symplectic) 2-form, respectfully denoted $\tilde{\theta},\,\tilde{\omega},$  by using the identification $P(T^*M)\cong T^*(PM).$ We can then pull these forms back to $\Pi_{\infty}^{(1)}(T^*M),$ and these give us $\theta,\,\omega.$ This make $\Pi_{2}(T^*M)$ into a presymplectic 2-groupoid\footnote{See \cite{vais} for more on presymplectic manifolds.} When the symplectic groupoid $\Pi_1(T^*M)$ exists, its symplectic form is obtained by reduction, and under certain conditions which we will discuss, the geometric quantization of $\Pi_{2}(T^*M)$ descends to a geometric quantization of $\Pi_1(T^*M)$ (see \cite{bon}, \cite{ruif} for details and explicit descriptions of the forms).
 \\\\That $dS=\delta^*\theta$ follows from first paragraph of proposition 2 in \cite{bon}. That is, consider a morphism $X:T\,\square\to T^*M$ ($\square$ is a square in the plane), diagrammatically given by the following:
\vspace{0.05cm}\\
\begin{equation*}
  \begin{tikzpicture}
\begin{scope}[very thick, every node/.style={sloped,allow upside down}]
\path
    (-2, 1) node[] {$X=$}
    (-1, 0) node[] {$\bullet$}
    (1, 0) node[] {$\bullet$}
    (1, 2) node[] {$\bullet$}
    (-1, 2) node[] {$\bullet$}
    (0,1)       node[] {$\Big\Uparrow$}
    ;
  \draw (1,0)-- node {\midarrow} (-1,0);
  \draw (1,0)-- node {} (1,2);
  \draw (1,2)-- node {\midarrow} (-1,2);
  \draw (-1,2)-- node {} (-1,0);

\end{scope}
\end{tikzpicture}
  \end{equation*}$\,$
  \vspace{0.05cm}\\
We can think of this morphism as defining a path $\gamma_s$ from the bottom arrow $\gamma_0$ to the top arrow $\gamma_1$, as indicated. We have that 
\vspace{0.05cm}\\
\begin{equation}
    \int_{\square}X^*\Pi=\int_{[0,1]}\gamma^*\theta
\end{equation}$\,$
\vspace{0.05cm}\\
The one-form $\theta$ defines a connection on the associated multiplicative line bundle, which we will discuss shortly.
\begin{remark}\label{reduc}
Actually, one could probably do a bit better by considering the simplicial space which in degree $n$ is given by all vector bundle morphisms $(X_n,\eta_n):T|\Delta|^n\to T^*M;$ denote this space by $V^n.$ From the above discussion $V^1$ has a 1-form and a (weak) symplectic form,\footnote{Weak meaning that the induced map $TV^1\to T^*V^1$ is injective, but not surjective.} while the first part of \ref{action} defines a real-vauled function on $V^2,$ ie. $(X_2,\eta_2)\to \int_{|\Delta|^2}\eta_2\wedge dX_2.$ Therefore, we get a cochain in degree $2$ on $V^{\bullet},$ and its pullback to $\Pi_{\infty}(T^*M)$ is our cocycle/geometric quantiztion. Thus, this cochain on $V^{\bullet}$ defines a kind of ``universal" geometric quantization. This should be compared with ``quantization commutes with reduction", see \cite{stern}. Also see \cite{felder} for details on how to perform the symplectic reduction.
\end{remark}
\begin{remark}
Assigning to $\Pi$ the map $X\mapsto\int_{\Delta}X^*\Pi $ should be viewed as a manifestation of the inverse of the van Est map on higher groupoids. That is, given any Lie algebroid $A\to M$ and any degree $n$ algebroid-form $\alpha$ on $A\,,$ we get a map on composable $n$-tuples of $\Pi_{\infty}(A),$ given by
\vspace{0.05cm}\\
\begin{equation}
    X^{(n)}\mapsto \int_{\Delta^n}X^{(n)*}\alpha\,,
\end{equation}$\,$
\vspace{0.05cm}\\
We will denote the map that associates to $\alpha$ the function above by 
\vspace{0.05cm}\\
\begin{equation}\label{int}
   \int: \alpha\mapsto \int \alpha\,,
   \end{equation}$\,$
   \vspace{0.05cm}\\
That $\int\alpha$ is a cocycle follows from Stokes' theorem and the fact that, since $X^{(n)}$ is a Lie algebroid morphism,
\vspace{0.05cm}\\
\begin{equation}
dX^{(n)*}\alpha=X^{(n)*}d_{\text{CE}}\,\alpha\,,
\end{equation}$\,$
\vspace{0.05cm}\\
where $d_{\text{CE}}$ is the Chevalley-Eilenberg differential. The map \ref{int} should define an isomorphism between the cohomology of $A$ and the cohomology of $\Pi_{\infty}(A),$ and on nice cocycles one can describe the inverse at the level of cochains. From this point of view, de Rham's theorem should be a special case of the van Est isomorphism theorem for higher groupoids, ie. letting $A=TX$ for some manifold $X\,,$ we know that \ref{int} is indeed an isomorphism as this is precisely what de Rham's theorem states.
\end{remark}$\,$
\vspace{0.05cm}\\
\subsection{The Simplicial Set $\mathbf{\Pi_{\infty}(T^*M\oplus_{\Pi}\mathbb{C})}$}
Now that we've described the cocycle we can describe the associated multiplicative line bundle $\mathcal{L}\to\Pi_{\infty}(T^*M)\,.$ Before doing this we will describe a simplicial set which we should think of as being the ``$\infty$-groupoid" (or 2-groupoid) integrating the canonical central extension determined by the 2-cocycle $\Pi\,$\footnote{See \ref{appen}.}:
\vspace{0.05cm}\\
\begin{equation}
    0\to M\times\mathbb{C}\to T^*M\oplus_{\Pi}\mathbb{C}\to T^*M\to 0\,.
\end{equation}$\,$
\vspace{0.05cm}\\
We will denote this simplicial set by $\Pi_{\infty}(T^*M)\oplus_{\Pi}\mathbb{C}^*\,.$
We will describe the simplicial set in low degree, in higher degrees it is similar: in degree $n$ the set is given by 
\vspace{0.05cm}\\
\begin{equation}\label{simp}
    \Pi_{\infty}^{(n)}(T^*M)\times \mathbb{C}^{*n}\,.
    \end{equation}
$\,$
\vspace{0.05cm}\\
The target and source maps on arrows are given by
\vspace{0.05cm}\\
\begin{equation}
d_0(\gamma,\lambda)=d_0(\gamma)\,,\,d_1(\gamma,\lambda)=d_1(\gamma)\,.
\end{equation}$\,$
\vspace{0.05cm}\\
On composable pairs the face maps are given by
\vspace{0.05cm}\\
\begin{align}
    &d_0(X,\lambda_1,\lambda_2)=(d_0(X),\lambda_2)\,,\,d_1(X,\lambda_1,\lambda_2)=(d_1(X),\lambda_1\lambda_2e^{iS[X]})\,,
    \\& d_2(X,\lambda_1,\lambda_2)=(d_2(X),\lambda_1)\,.
\end{align}$\,$
\vspace{0.05cm}\\
The face maps are most important in degree one and two. On composable triples the face maps are
\vspace{0.05cm}\\
\begin{align}
    &d_0(X^{(3)},\lambda_1,\lambda_2,\lambda_3)=(d_0(X^{(3)}),\lambda_2,\lambda_3)\,,\\
    &\,d_1(X^{(3)},\lambda_1,\lambda_2,\lambda_3)=(d_1(X^{(3)}),\lambda_1\lambda_2e^{iS[d_3(X^{(3)})]},\lambda_3)\,,\\
    &d_2(X^{(3)},\lambda_1,\lambda_2,\lambda_3)=(d_2(X^{(3)}),\lambda_1,\lambda_2\lambda_3 e^{iS[d_0(X^{(3)})]})\,,\\
    &d_3(X^{(3)},\lambda_1,\lambda_2,\lambda_3)=(d_3(X^{(3)}),\lambda_1,\lambda_2)
\end{align}$\,$
\vspace{0.05cm}\\
\subsection{The Multiplicative Line Bundle $\mathbf{\mathcal{L}\to \Pi_{\infty}(T^*M)}$}
Now we get a new simplicial set by replacing the $\mathbb{C}^*$ in \ref{simp} with $\mathbb{C}\,.$  That is, in degree $n$ the set is given by 
\vspace{0.05cm}\\
\begin{equation}
    \Pi_{\infty}^{(n)}(T^*M)\times \mathbb{C}^{n}\,.
    \end{equation}$\,$
    \vspace{0.05cm}\\
and the simplicial maps are otherwise the same as in the previous section. There is a natural map 
\vspace{0.05cm}\\
\begin{equation}
    \Pi_{\infty}^{(n)}(T^*M)\times \mathbb{C}^{n}\to \Pi_{\infty}^{(n)}(T^*M)
    \end{equation}$\,$
    \vspace{0.05cm}\\
in each degree. We will think of this as a multiplicative line bundle and denote it by $\mathcal{L}\to \Pi_{\infty}(T^*M)\,.$
\subsection{The Twisted Convolution Algebra}
Now that we have defined the multiplicative line bundle we are ready to define the twisted convolution algebra. First, a definition:
\vspace{0.05cm}\\
\begin{definition}
A section of $\mathcal{L}\to \Pi_{\infty}(T^*M)$ is defined to be a section of $\mathcal{L}^{(1)}\to \Pi_{\infty}^{(1)}(T^*M)\,.$ We may call these sections of $\mathcal{L}$ for brevity. These can be identified with maps 
\vspace{0.05cm}\\
\begin{equation}
    \text{Hom}(T\Delta^1,T^*M)\to \mathbb{C}\,.
\end{equation}
\end{definition}$\,$
\vspace{0.05cm}\\
Now we will define the twisted convolution algebra similarly to how it is done on Lie groupoids.
\vspace{0.05cm}\\
\begin{definition}
Let $s_1\,,s_2$ be sections of $\mathcal{L}\,.$ Their convolution is another section of $\mathcal{L}\,,$ defined by\\
\begin{equation}\label{con}
    (s_1 \ast s_2)(\gamma)=\int_{d_1(X)=\gamma}s_1(d_0(X))s_2(d_2(X))\,e^{\frac{i}{\hbar}S[X]}\,dX\,,
\end{equation}$\,$
\vspace{0.05cm}\\
where the measure is normalized so that $1\star 1=1$ (where $1$ is the constant section $\gamma\mapsto (\gamma,1))\,.$
\end{definition}$\,$
\vspace{0.05cm}\\
\begin{remark}
If we take $s_1=\lambda_1\delta_{\gamma_1}\,,s_2=\lambda_2\delta_{\gamma_2}\,,$ then we can formally compute 
\vspace{0.05cm}\\
\begin{equation}
(\lambda_1\delta_{\gamma_1}\ast\lambda_2\delta_{\gamma_2})(\gamma)=\lambda_1\lambda_2\int_{\begin{subarray}{l}d_0(X)=\gamma_2,\\d_1(X)=\gamma, d_2(X)=\gamma_1\end{subarray}}e^{\frac{i}{\hbar}S[X]}\,dX\,.
\end{equation}$\,$
\vspace{0.05cm}\\
This is formally similar to the product in the Fukaya category if $(M,\Pi)$ is symplectic and if $\gamma_1,\gamma_2,\gamma$ are Lagrangian submanifolds (see \cite{zas}).
\end{remark}$\,$
\vspace{0.05cm}\\
\section{The Product}\label{starproduct}
Now that we have defined the twisted convolution algebra we can use it define a product on functions in $C^{\infty}(M).$ In order to do this, we need to associate sections of $\mathcal{L}$ to functions in $C^{\infty}(M)\,.$ To do this, we choose a generic point on $|\Delta^1|$ (the interval), ie. any point that isn't an endpoint. We will take this generic point to be $t=1/2\,.$\\
\vspace{0.05cm}\\
Let $\pi:T^*M\to M$ denote the projection. Given $f\in C^{\infty}(M)\,,$ we get a section $\hat{f}$ of $\mathcal{L}$ by defining
\vspace{0.05cm}\\
\begin{equation}
    \hat{f}(\gamma):=f(\pi(\gamma(1/2)))\,,
\end{equation}
\vspace{0.05cm}\\
and this identifies $C^{\infty}(M)$ with a subspace of the sections of $\mathcal{L}\to \Pi_{\infty}(T^*M)\,.$ Notice that for $m\in M,$ thought of as the constant path, $\hat{f}(m)=f(m).$ Now we are ready to define a new product on $C^{\infty}(M)\,.
\vspace{0.05cm}\\
$ \begin{definition}We define a product on $C^{\infty}(M)$ by
\vspace{0.05cm}\\
\begin{equation}\label{starr}
    (f\star g)(m):=(\hat{f}\ast \hat{g})(m)\,,
\end{equation}
\vspace{0.05cm}\\
where $m\in M$ and where on the right side of the above equation we use the identification of $C^{\infty}$ with sections of $\mathcal{L}\,.$
\end{definition}$\,$
\vspace{0.05cm}\\
Explicitly, we have that 
\vspace{0.05cm}\\
\begin{align}\label{mystar}
    (f\star g)(m)&=\int_{X\text{ composes to }m}\hat{f}(\gamma_2)\,\hat{g}(\gamma_1)\,e^{\frac{i}{\hbar}S[X]}\,dX\\&=\int_{X\text{ composes to }m}f(\pi(\gamma_2(1/2)))\,g(\pi(\gamma_1(1/2)))\,e^{\frac{i}{\hbar}S[X]}\,dX\,,
\end{align}$\,$
\vspace{0.05cm}\\
where the measure is normalized so that $1\star 1=1\,.$
This is the same product as in equation 3.45 in \cite{bon}.
\vspace{0.05cm}\\
\section{Finite Dimensional Phase Spaces}
Soon we will discuss a way of reducing the product $\star$ to a finite dimensional integral (for gauge invariant observables, which we will discuss shortly).
Before doing so, we are going to rewrite the integral. In the followig diagram, consider the composable pair $X$ on the left, for which $\gamma_2\circ\gamma_1\Rightarrow \text{Id}_m\,,$ and with the marked points for each arrow at $t=1/2\,.$ We can choose a diffeomorphism of the disk $h:D^2\to D^2$ which fixes the marked points in such a way that the composition $X\circ h$ gives the composable pair $\tilde{X}$ on the right, for which
$\tilde{\gamma_2}\circ\tilde{\gamma_1}\Rightarrow\tilde{\gamma}_3\,.$
\vspace{0.05cm}\\
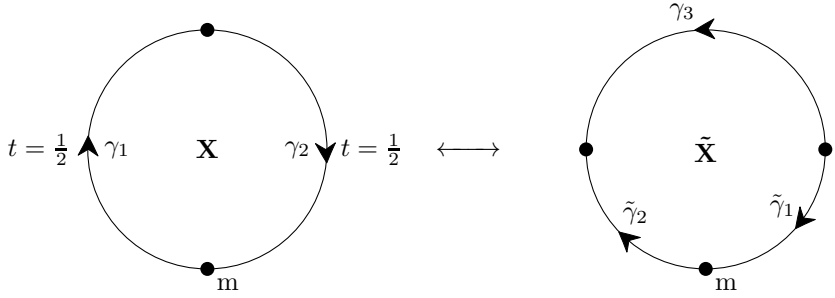
\begin{adjustwidth*}{}{}
\begin{tikzpicture}[scale=5.3,cap=round,>=latex]
\def\Radius{.3cm}
\draw (0cm,0cm) circle[radius=\Radius];

\begin{scope}[
    -{Stealth[round, length=8pt, width=8pt, bend]},
    shorten >=4pt,
    very thin,
  ]
    \draw (-5:\Radius) arc(185:187:1);
    \draw (175:\Radius) arc(361:363:1);
  \end{scope}

  \fill[radius=0.5pt]
  (-90:\Radius) circle[] node[above left] {}
    (90:\Radius) circle[] node[above left] {}

    ;

      \path
    (0:\Radius) node[right] {$\,t=\frac{1}{2}$}
    (180:\Radius) node[ left] {$t=\frac{1}{2}\,\,$}
    (-90:\Radius) node[below right] {m}
    (0.55,0.05) node[below right] {$\xleftrightarrow[]{\;\;\;\;\;\;\;\;}$}
    (0:\Radius) node[left] {$\gamma_2\;$}
    (180:\Radius) node[right] {$\;\gamma_1\;$}
    (0:0) node[] {$\mathbf{X}$}
    
  ;

\draw (1.25cm,0cm) circle[radius=\Radius];
\begin{scope}[
    -{Stealth[round, length=8pt, width=8pt, bend]},
    shorten >=4pt,
    very thin,
  ]
    \draw (0.2298+1.25, -0.1928) arc(143:145:1);
    \draw (-0.02+1.25, 0.3) arc(95:97:1);
    \draw (-0.212132+1.25, -0.212132) arc(50:52:1);
  \end{scope}

  \fill[radius=0.5pt]
    (1.25, -\Radius) circle[] node[above left] {}
     (0.95, 0) circle[] node[above left] {}
     (1.55, 0) circle[] node[above left] {}

    ;

      \path
    (-45:\Radius+1.2) node[below right] {}
    (-135:\Radius+1.25) node[above left] {}
    (1.25, -\Radius) node[below right] {m}
    (0.2298+1.25, -0.1928) node[above] {$\tilde{\gamma}_1\;\;\;\;$}
    (-0.05+1.25, 0.3) node[above] {$\tilde{\gamma}_3\;$}
    (-0.212132+1.25, -0.212132) node[above] {$\;\;\;\;\;\tilde{\gamma}_2\;$}
    (1.25,0) node[] {$\mathbf{\tilde{X}}$}
    
  ;

\end{tikzpicture}
\end{adjustwidth*}$\,$
\vspace{0.05cm}\\
Essentially, we are just thinking of $X$ as being a composable pair in a different way; this is similar to the diffeomorphism of $G^{(2)}$ given by $(g_2,g_1)\mapsto (g_2^{-1},g_2g_1),$ for a Lie groupoid $G.$ Therefore, rather than integrating over all composable pairs on the left, we can integrate over all composable pairs on the right, ie. we integrate over all composable pairs $\tilde{X}$ such that $\pi(\tilde{\gamma}_1(1))=\pi(\tilde{\gamma}_2(0))=m\,.$ Let $s(\gamma)=\pi(\gamma(0))\,,t(\gamma)=\pi(\gamma(1))\,.$ Then, we can rewrite the product as 
\vspace{0.05cm}\\
\begin{align}\label{starry}
    (f\star g)(m)=\int_{\tilde{X}}t^*f(\tilde{\gamma}_2)\,s^*g(\tilde{\gamma}_1)\,e^{\frac{i}{\hbar}S[\tilde{X}]}\,d\tilde{X}\,.
\end{align}    
\vspace{0.05cm}\\
 The advantage of this is that if 
 \begin{equation}\label{condo}
     \int_{S^2}\tilde{X}^*\Pi=0\,,
\end{equation}
\vspace{0.05cm}\\
for all algebroid morphisms $\tilde{X}:TS^2\to T^*M\,,$ then the integrand of $f\star g$ is now written in terms of quantities which make sense on the source simply connected integration $\Pi_1(T^*M)\,.$\\
 \vspace{0.05cm}\\The sense in which $\tilde{X}\mapsto e^{\frac{i}{\hbar}S[\tilde{X}]}$ descends to the groupoid is the following: while the cocycle itself doesn't descend as an $S^1$-valued map on composable pairs, it does descend as a multiplicative line bundle (we will discuss this in more detail shortly). This line bundle has a trivialization which induces an $S^1$-valued cocycle $e^{\frac{i}{\hbar}S_0[\cdot,\cdot]}$ on composable arrows in $G\,.$ This cocycle can always be chosen to have the following properties:  
 \vspace{0.05cm}\\
 \begin{enumerate}
       \item $VE(S_0)=\Pi\,,$
       \item $S_0[g_2,g_1]=-S_0[g_1^{-1},g_2^{-1}]\,.$
\end{enumerate}$\,$
\vspace{0.05cm}\\
The first condition says that the van Est map (see \ref{appen}) applied to $S_0$ gives the algebroid cocycle $\Pi\,,$ while the second condition is an antisymmetry condition, which is analogous to the fact that  $X\mapsto\int_{|\Delta|^2}X^*\Pi$ changes sign when reversing the orientation of $|\Delta^2|\,.$ This should be thought of as being obtained by symplectic reduction, see remark \ref{reduc}.   
\vspace{0.05cm}\\
\section{$\mathbf{C^*}$-Algebras of Poisson Manifolds}
\noindent Here we are going to explain a method of assigning a $C^*$-algebra to a Poisson manifold, which will be done by interpreting \ref{starry} as an integral over $\Pi_1(T^*M)$ (the source simply connected integration). For now we will assume that $(M,\Pi)$ is integrable, and that
\vspace{0.05cm}\\
\begin{equation}\label{conditionnn}
  \int_{S^2}X^*\Pi=0   
 \end{equation}
 \vspace{0.05cm}\\
 for all algebroid morphisms $X:TS^2\to T^*M\,.$ We will discuss the more general case later. For now we need a cocycle and a measure. To see an example fully carried out, see \ref{fund}.
\vspace{0.05cm}\\
\subsection{The Cocycle}\label{cocycle}$\,$
 \vspace{0.05cm}\\
By the van Est isomorphism theorem (see \ref{VE}), the condition \ref{conditionnn} is necessary and sufficient for the existence of a degree 2 cocycle $S$ on $\Pi_1(T^*M)$ such that 
 \vspace{0.05cm}\\
 \begin{equation}
     VE(S)=\Pi\,,
 \end{equation}
  \vspace{0.05cm}\\
 where $VE$ is the van Est map, applied at the level of cochains. Now given such a cocycle, we get a new cocycle $\tilde{S}$ by defining
  \vspace{0.05cm}\\
 \begin{equation}
     \tilde{S}(g_2,g_1)=\frac{S(g_2,g_1)-S(g_1^{-1},g_2^{-1})}{2}\,.
 \end{equation}
  \vspace{0.05cm}\\
 We still have that $VE(\tilde{S})=\Pi\,,$ and we further have the antisymmetry property
  \vspace{0.05cm}\\
 \begin{equation}
     \tilde{S}(g_2,g_1)=-\tilde{S}(g_1^{-1},g_2^{-1})\,,
 \end{equation}
  \vspace{0.05cm}\\
 therefore we may always assume the cocycle is antisymmetric. This is similar to a comment made by Kontsevich on page 3 of \cite{kontsevich}: the symmetric part of the first order differential operator in the star product can always be killed by a gauge transformation; at least if our groupoid is source simply connected any symmetric 2-cocycle is cohomologically trivial, due to the the van Est isomorphism theorem. We will come back to this in example \ref{good}.
 \\\\The benefit of having an antisymmetric cocycle is that the involution of the $C^*$-algebra we are going to construct will just be complex conjugation. We then get an $S^1$-valued cocycle by forming  $e^{\frac{i}{\hbar}S}\,,$ and this is the analogue of $e^{\frac{i}{\hbar}\int_{\Delta}X^*\Pi}\,.$ 
 \subsection{The Measure}$\,$
  \vspace{0.05cm}\\
 Now that we have a cocycle, we need to discuss the analogue of the measure $d\tilde{X}$ in \ref{starry}. Given an $m\in M\,,$ this will be a measure on the on the space of composable arrows ending and beginning at $m\,,$ ie.
  \vspace{0.05cm}\\
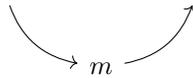
 \begin{equation*}
 \begin{tikzcd}
{} \arrow[rd, bend right] &                          & {} \\
                          & m \arrow[ru, bend right] &   
\end{tikzcd}
 \end{equation*}
  \vspace{0.05cm}\\
 To find this measure, choose a top form on $T^*M\,.$ This induces left and right Haar measures $d\mu$ on $\Pi_1(T^*M)$ by left and right translation\footnote{As we will see soon, sometimes there are better ways of choosing invariant measures.}. This gives a measure on each source and target fiber, and this induces a measure on $s^{-1}(m)\times t^{-1}(m)\,,$ for each $m\in M\,.$ We will denote this measure by $d\mu_m\otimes d\mu_m\,.$ 
 \vspace{0.5cm}\\Let us remark that, in the case that each isotropy group is \textit{unimodular}, the Haar measure can be chosen in a unique way up to a Casimir function, ie. a function constant on symplectic leaves. To see this, pick an $m\in M\,.$ Given $g\in s^{-1}(m)\,,$ the isotropy group at $t(g)$ acts freely and transitively on the fiber of 
 \vspace{0.05cm}\\
 \begin{equation}
 t\vert_{s^{-1}(m)}:s^{-1}(m)\to M
  \end{equation}
   \vspace{0.05cm}\\
 over $t(g)\,,$ and therefore a Haar measure on the isotropy group over $t(g)$ induces a measure on this fiber. Thus, smoothly choosing a Haar measure for each isotropy group over the symplectic leaf $\mathcal{L}_m$ containing $m$ induces a measure along the fibers of $t\vert_{s^{-1}(m)}\,.$ 
 \vspace{0.5cm}\\Furthermore, the symplectic leaf has a natural volume form associated to its symplectic form, denoted $\omega^n\,,$ where $n=(\text{dim }\mathcal{L}_m)/2\,.$ Together, the measures along the symplectic leaf and the fibers induce a measure on the total space, ie. a measure on $s^{-1}(m)\,.$ Doing this for each $m\in M$ determines a right Haar measure on $G\,.$ Similarly, we get a left Haar measure.
 \vspace{0.5cm}\\Now, the Haar measure described in the previous paragraph is only defined up to a function on $M\,.$ However, given any two points in a symplectic leaf and an arrow between them, we get an isomorphism between the isotropy groups, and any two such isomorphisms differ by an inner automorphism. Since the istotropy groups are assumed to be unimodular, their Haar measures are bi-invariant and therefore invariant under inner automorphisms. Therefore we get the following:
  \vspace{0.05cm}\\
 \begin{lemma}
Let $(G\rightrightarrows M, \omega)$ be a symplectic groupoid over $(M,\Pi)$ such that all isotropy groups are unimodular. Then up to a Casimir function there is a canonical choice of Haar measure.
 \end{lemma}
 Of course, a compact Lie group has a canonical Haar measure, and since compact Lie groups are unimodular we find that:
 \begin{lemma}
Let $(G\rightrightarrows M, \omega)$ be a symplectic groupoid over $(M,\Pi)$ such that all isotropy groups are compact. Then there is a canonical choice of Haar measure.
 \end{lemma}$\,$\\
 \begin{exmp}
Let $(M,\omega)$ be a symplectic manifold. Its isotropy groups are $\{e\}$ and therefore they are compact and come with a unique normalized measure. We have that $s^*\omega,\,t^*\omega$ are left and right Haar measures on $\Pi_1(M),$ respectively. 
 \end{exmp}$\,$
  \vspace{0.05cm}\\
 In the next subsection we will (in particular) discuss normalization of the measures $d\mu_m\,.$
  \vspace{0.05cm}\\
 \subsection{The Product}$\,$
  \vspace{0.05cm}\\
Let $(M,\Pi)$ be a Poisson manifold and $G\rightrightarrows M$ a Lie groupoid integrating it. Then given a cocycle $S$ satisfying \ref{condo}, together with a measure $d\mu_m$ on $s^{-1}(m)\times t^{-1}(m)$ for each $m\in M\,,$ as discussed in the previous subsection, we \textit{formally} define a product by\footnote{Formally, this is the same as \ref{starry}.} \\
 \vspace{0.05cm}\\
\begin{equation}\label{guess1}
    \boxed{(f\star g)(m)=\frac{1}{C_m}\int_{s^{-1}(m)\times t^{-1}(m)} d\mu_m\otimes d\mu_m \,(t^*f)(s^*g) e^{\frac{i}{\hbar}S}}\,,
\end{equation}
 \vspace{0.1cm}\\\\
where $C_m$ is a normalization parameter and is defined so that $1\star 1=1$ (where $1\star 1$ is interpreted as an improper integral). It is given by 
 \vspace{0.05cm}\\
\begin{equation}
    C_m=\int_{s^{-1}(m)\times t^{-1}(m)} d\mu_m\otimes d\mu_m \,e^{\frac{i}{\hbar}S}\,.
\end{equation}
 \vspace{0.05cm}\\
We use the term \textit{formally} because there may be symmetries causing some divergence issues for \ref{guess1}. In all of the examples we see the cocycle is invariant under the actions of the identity components of the isotropy groups — as is $\Pi.$ There are two ways of dealing with this: one is to let the invariant measures on the isotropy groups be only finitely additive so that we can use (finitely additive) probability measures,\footnote{See \cite{Runde}.} and another is to quotient out the domain of integration (ie. the phase space) by the symmetries and integrate  over the quotient. We will further discuss this with the following example:
 \vspace{0.05cm}\\
\begin{exmp}\label{basic}
The simplest Poisson manifold is the $0$-Poisson structure on a manifold $M\,.$ The groupoid is the vector bundle $T^*M\rightrightarrows M$ and the cocycle is zero. In this case the symplectic leaves are points and the isotropy is $\mathbb{R}^{n},$ where $n=\text{dim M}\,.$ For the Haar measure, one option is the Lebesgue measure $d^n x\,.$ Equation \ref{guess1} reduces to
 \vspace{0.05cm}\\
\begin{align}
    (f\star g)(m)&=\frac{1}{C_m}\int_{\mathbb{R}^{n}\times \mathbb{R}^{n}} d^{n}x\,d^n{y} \,f(m)g(m)\\
    &=\frac{f(m)g(m)}{C_m}\int_{\mathbb{R}^{n}\times \mathbb{R}^{n}} d^{n}x\,d^n{y}\,,
\end{align}
 \vspace{0.05cm}\\
and since formally $C_m=\int_{\mathbb{R}^{n}\times \mathbb{R}^{n}} d^{n}x\,d^n{y}\,,$ this just gives $f(m)g(m)\,.$ However since the normalization parameter $C_m=\infty\,,$ this tells us that we shouldn't be using the Lebesgue measure. 
\paragraph{Method 1:}$\,$\\\vspace{0.05cm}\\ Notice that, if we instead integrate over a bounded domain and take a limit, everything will work fine. Therefore, we will interpret these equations as telling us that we should be using an invariant probability measure on $\mathbb{R}^n$ instead of the Lebesgue measure. Since $\mathbb{R}^n\rightrightarrows *$ is amenable, it has a finitely additive, invariant probability measure (or invariant mean) $d\mu\,,$ and redoing the computation with this measure we get that $C_m=1$ and 
 \vspace{0.05cm}\\
\begin{align}
    (f\star g)(m)&=\int_{\mathbb{R}^{n}\times \mathbb{R}^{n}} d\mu(x)\,d\mu(y)\,f(m)g(m)=f(m)g(m)\,.
\end{align}
 \vspace{0.05cm}\\
\paragraph{Method 2:}$\,$\\\vspace{0.05cm}\\ The second method of dealing with the infinty is to further reduce the phase space of the integral. In this case, since the observables and the cocycle are constant along the isotropy, we could quotient out by the entire isotropy group and thus we only need to integrate over the point $\{m\}\,,$ which has a unique probability measure $\mu\,.$ Doing this we get 
 \vspace{0.05cm}\\
\begin{align}
    (f\star g)(m)&=\int_{\{m\}}\,d\mu\,f(m)g(m)=f(m)g(m)\,.
\end{align}
 \vspace{0.05cm}\\
Therefore, both methods give the expected $f\star g=fg\,.$
\end{exmp}$\,$
 \vspace{0.05cm}\\
\subsection{Dealing with Infinities}$\,$
 \vspace{0.05cm}\\
Let's expound on the two methods of dealing with infinities discussed in the previous subsection, which can occur when the cocycle is invariant under a subgroup of the isotropy.
\paragraph{Method 1: Finitely additive probability measures}$\,$\\
\vspace{0.05cm}\\On a Lie group there is a unique left (and right) invariant measure, up to a constant\footnote{Usually this is formulated as an invariant mean.}; this is called the Haar measure. However, if one relaxes the notion of countable additivity to finite additivity this uniqueness no longer necessarily holds. An amenable group is a group with a finitely additive, left (or right) invariant probability measure with respect to the sigma algebra associated with the Haar measure, eg. abelian groups, compact groups, nilpotent groups (see \cite{Runde}). As seen in example \ref{basic}, sometimes these finitely additive (left or right) invariant probability measures can be used rather than the Haar measures to deal with infinities.
\paragraph{Method 2: Quotient Out by Gauge Symmetries}$\,$\\
\vspace{0.05cm}\\
Let $m\in M\,.$ The integrand of $(f\star g)(m)$ is given by 
\vspace{0.05cm}\\
\begin{equation*}
    t^*f(\gamma_2)s^*g(\gamma_1)\ e^{\frac{i}{\hbar}S[\gamma_2,\gamma_1]}\,.
    \end{equation*}
    \vspace{0.05cm}\\
Now, because the isotropy $G_m$ of $G\rightrightarrows M$ over $m$ acts on both $t^{-1}(m)\,, s^{-1}(m)\,,$ $G_m\times G_m$ acts on $t^{-1}(m)\times s^{-1}(m)\,,$ and therefore it acts on the integrand, via\footnote{Here $(k_1,k_2)\in G_m\times G_m$ and $(\gamma_2,\gamma_1)\in t^{-1}(m)\times s^{-1}(m).$}
\vspace{0.05cm}\\
\begin{align}
    &\big[(k_2,k_1)\cdot [t^*f][s^*g]\ e^{\frac{i}{\hbar}S}\big](\gamma_2,\gamma_1)
=t^*f(\gamma_2\cdot k_2)s^*g(k_1\cdot\gamma_1)\ e^{\frac{i}{\hbar}S[\gamma_2\cdot k_2,k_1\cdot\gamma_1]}\,. 
\end{align}
\vspace{0.05cm}\\
Of course, $G_m\times G_m$ fixes $t^*f\,s^*g\,,$ but in general it won't fix $e^{\frac{i}{\hbar}S}\,.$
However, suppose there is a subgroup $H_m\subset G_m$ such that $H_m\times H_m$ fixes  $e^{\frac{i}{\hbar}S}\,.$ Then, we could try to integrate over 
\vspace{0.05cm}\\
\begin{equation}
    t^{-1}(m)\times s^{-1}(m)/(H_m\times H_m)\cong  t^{-1}(m)/H_m\times s^{-1}(m)/H_m
\end{equation}
\vspace{0.05cm}\\
instead. Note that, after taking closures, we may assume that $H_m\subset G_m$ is closed. Recall that the map
\vspace{0.05cm}\\
\begin{equation}
    t\vert_{s^{-1}(m)}:s^{-1}(m)\to \mathcal{L}_m
\end{equation}
\vspace{0.05cm}\\
is a principal $G_m$-bundle, where $\mathcal{L}_m$ is the symplectic leaf through $m\in M\,.$ Taking the fiberwise quotient by $H_m\,,$ we get a fiber bundle with fiber $G_m/H_m:$
\vspace{0.05cm}\\
\begin{equation}
    t\vert_{s^{-1}(m)/H_m}:s^{-1}(m)/H_m\to \mathcal{L}_m\,.
\end{equation}
\vspace{0.05cm}\\
Of course, $\mathcal{L}_m$ still has its natural measure associated with the symplectic form. Therefore, in order to get a measure on $s^{-1}(m)/H_m$ we just need a measure on $G_m/H_m\,.$ We have the following theorem:
\vspace{0.05cm}\\
\begin{theorem}(see \cite{Gorb}, chapter 1)
Let $G$ be a Lie group and $H\subset G$ a closed subgroup.  Then there is a $G$-invariant measure $\mu_{G/H}$ on $G/H$ if and only if $\Delta_G\vert_H=\Delta_H\,,$ where $\Delta_G,\,\Delta_H$ are the respective modular functions. In this case, $\mu_{G/H}$ has the following property: for any $f\in L^1(G)\,,$ we have that
\vspace{0.05cm}\\
\begin{equation}
    \int_{G/H}\Big(\int_H f(gh)d\mu_H\Big)d\mu_{G/H}=\int_G f(g)\,d\mu_G\,,
\end{equation}
\vspace{0.05cm}\\
ie. integrating along the fiber $H$ and then integrating along the base $G/H$ is the same as integrating over $G\,.$
\end{theorem}$\,$
\vspace{0.05cm}\\
The above theorem implies, in particular, that if both $H,\,G$ are unimodular then such a measure exists on $G/H\,.$ This is always the case if $G$ is compact, abelian or nilpotent, for example.
\vspace{0,5cm}\\In example \ref{basic}, the cocycle was invariant under the entire isotropy group, and thus 
\vspace{0.05cm}\\
\begin{equation}
    t^{-1}(m)\times s^{-1}(m)/(H_m\times H_m)\cong \{m\}\,.
    \end{equation}
    \vspace{0.05cm}\\
    \subsection{The $\mathbf{C^*}$-Norm and Involution}$\,$\\
    Finally, assuming we successfully dealt with the infinities to get convergent integrals, the involution is given by complex conjugation (this is due to the antisymmetry condition on the cocycle), and the $C^*$-norm is given by the operator norm (this is very similar to the reduced $C^*$-algebra norm on a groupoid, see \cite{bu}, \cite{eli}).
    \subsection{On Associativity of the Product}$\,$\\
    The product \ref{guess1} is associative exactly when it is associative on delta functions, ie. associativity is equivalent to\\
    \begin{equation}\label{delta}
        (\delta_a\star\delta_b)\star\delta_c=  \delta_a\star(\delta_b\star\delta_c)\,,
    \end{equation}\\
for all $a,b,c\in M.$
We can think of this in the following way: for all $x\in M,$ it is necessary and sufficient that\footnote{This is assuming we have dealt with any infinities.}
\begin{align}\label{ass}
   & \int_{\Omega_1}e^{\frac{i}{\hbar}[S(\gamma_{a},\,\gamma_{b})+S(\gamma_{x'},\,\gamma_{c})]}\;d\gamma_a\,d\gamma_b\,d\gamma_c\,d\gamma_{x'}\,\omega^n_{x'}\nonumber
   \\&=\int_{\Omega_2}e^{\frac{i}{\hbar}[S(\gamma_{a},\,\gamma_{x'})+S(\gamma_{b},\,\gamma_{c})]}\;d\gamma_a\,d\gamma_b\,d\gamma_c\,d\gamma_{x'}\,\omega^n_{x'}\,,
\end{align}$\,$\\
where $n$ is half the dimension of the symplectic leaf through $x\,.$ To be clear, $a, b ,c ,x$ are fixed.
The domains of integration $\Omega_1,\,\Omega_2$ are given by all diagrams of the following form, where the $\gamma$'s are arrows in the groupoid, and the $d\gamma$'s are the induced measures obtained from \ref{guess1}:\\
\begin{equation}
  \Omega_1:\;\;  \begin{tikzcd}
                            &                           & a                        \\
                            & x' \arrow[ru, "\gamma_a"] &                          \\
x \arrow[ru, "\gamma_{x'}"] &                           & b \arrow[lu, "\gamma_b"] \\
                            & c \arrow[lu, "\gamma_c"]  &                         
\end{tikzcd}
\;\;\;\;\;\;\;\;\;\;\;\;\;\;\;
\Omega_2:\;\;\begin{tikzcd}
                         & a                                                   &                          \\
x \arrow[ru, "\gamma_a"] &                                                     & b                        \\
                         & x' \arrow[lu, "\gamma_{x'}"] \arrow[ru, "\gamma_b"] &                          \\
                         &                                                     & c \arrow[lu, "\gamma_c"]
\end{tikzcd}
\end{equation}$\,$\\\\
In example \ref{symp} we will use this to show associativity of the Moyal product.
\section{Relationship Between Groupoid and Poisson Sigma Models}$\,$
    \vspace{0.05cm}\\
Up until now we haven't much discussed the sense in which the groupoid-valued sigma model computes correlation functions of gauge invariant observables. In \cite{catt} the authors consider two morphisms $X_1,X_2:T|\Delta^2|\to T^*M$ to be gauge equivalent if $\partial X_1=\partial X_2:T\partial|\Delta^2|\to M\,,$  ie. two morphisms are equivalent if they agree on the boundary. For a thorough discussion of some gauge groups relevant to the Poisson sigma model, see \cite{bojo}, \cite{felder}, \cite{zab}, \cite{catt}. Here we will sketch the idea of our gauge equivalences. Basically, the relationship goes: under symplectic reduction,\\
\begin{align}
    &\text{Geometric quantization of }\Pi_2(T^*M)\leadsto\text{Geometric quantization of }\Pi_1(T^*M)\,,\\
    &f\star g \text{ with respect to }\Pi_2(T^*M) \leadsto f\star g \text{ with respect to }\Pi_1(T^*M)\,,
\end{align}\\
Recall that the reduced phase space of $\Pi_{\infty}^{(1)}(T^*M)$ is $\Pi_{1}^{(1)}(T^*M)$ (see \cite{felder}), so we are essentially performing reduction on the quantization we constructed on $\Pi_2(T^*M).$ See section \ref{idea} and remark \ref{reduc}.\\
\vspace{0.05cm}\\We want to work with a slightly coarser gauge equivalence. Rather than considering two morphisms $X_1,X_2:T|\Delta^2|\to T^*M$ to be equivalent only if they agree on the boundary, we consider them to be equivalent if the three boundary components are Lie algebroid homotopy equivalent, with fixed endpoints. That is $d_iX_1\sim d_i X_2$ as maps $Td_i|\Delta^2|\to T^*M,\,i=0,1,2,$ where the endpoints of $d_i|\Delta^2|\,,i=0,1,2$ are fixed throughout the homotopy. Equivalently, we consider $X_1,X_2:T|\Delta^2|\to T^*M$ to be equivalent if their images in $\Pi_1^{(2)}(T^*M)$ under the map $\Pi_{\infty}(T^*M)\to \Pi_{1}(T^*M)$ are equal.\\
\vspace{0.05cm}\\Now let's consider the multiplicative line bundle on $\Pi_{\infty}(T^*M)$ induced by the cocycle $e^{\frac{i}{\hbar}S[\tilde{X}]}\,.$ As shown in \cite{catt}, under the assumption that 
\vspace{0.05cm}\\
\begin{equation}\label{condition}
    \int_{S^2}X^*\Pi=2\pi\hbar\mathbb{Z},
\end{equation}
\vspace{0.05cm}\\
for all morphisms $X:TS^2\to T^*M\,,$ this multiplicative line bundle desdends to $\Pi_1(T^*M)\,.$ This can be seen by the following construction (see \cite{catt} for a similar construction): consider the space $(\Pi_{\infty}^{(1)}(T^*M)\times\mathbb{C})/\sim,$ where \\
\begin{equation}
    (\lambda,\gamma)\sim (\lambda e^{\frac{i}{\hbar}S[X]},\gamma')
    \end{equation}$\,$\\
    
if $X:TD\to T^*M$ is a morphism from the disk with two marks points, and with boundaries $\gamma,\,\gamma',$ as in the the following picture:
\vspace{0.05cm}\\
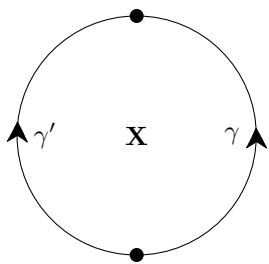
\begin{center}
\begin{tikzpicture}[scale=5.3,cap=round,>=latex]
\def\Radius{.3cm}
\draw (0cm,0cm) circle[radius=\Radius];

\begin{scope}[
    -{Stealth[round, length=8pt, width=8pt, bend]},
    shorten >=4pt,
    very thin,
  ]
    \draw (2:\Radius) arc(3:5:1);
    \draw (175:\Radius) arc(361:363:1);
  \end{scope}

  \fill[radius=0.5pt]
  (-90:\Radius) circle[] node[above left] {}
    (90:\Radius) circle[] node[above left] {}

    ;

      \path
  
    (0:\Radius) node[left] {$\gamma\;$}
    (180:\Radius) node[right] {$\;\gamma'\;$}
    (0:0) node[] {$\mathbf{X}$}
    
  ;

\end{tikzpicture}
\end{center}$\,$
\vspace{0.05cm}\\
Now by the assumption \ref{condition}, no two points in the vector space over the same $\gamma$ get identified, and this is why it descends to a line bundle. We will denote points in the equivalence classes by $[(\lambda,\gamma)]\,.$ The projection is given by
\vspace{0.05cm}\\
\begin{equation}
    [(\lambda,\gamma)]\mapsto [\gamma]\,.
\end{equation}
The multiplication is defined by
\begin{equation}
    [(\lambda,\gamma_2)]\cdot[(\beta,\gamma_1)]=[(\lambda\beta e^{\frac{i}{\hbar}S[X]},\gamma_3)],
\end{equation}
\vspace{0.05cm}\\
where $[\gamma_3]=[\gamma_2]\cdot[\gamma_1]\,,$ and $X:T|\Delta^2|\to T^*M$ is any morphism composing $\gamma_1,\gamma_2$ to give $\gamma_3.$ This product is well-defined, and this defines a multiplicative line bundle on $\Pi_1(T^*M)\,.$\footnote{The existence of the multiplicative line bundle is proved by cohomologically means in \cite{Lackman}, and a different proof appears in \cite{zhuc}.}\\
\vspace{0.05cm}\\Now one subtle point is that, looking at \ref{starry}, we are implicitly using the constant section $s(\gamma)=1$ to define the observables. The problem is that the integrand is not actually gauge invariant, which is due to the fact that the constant section is not gauge invariant, ie. it doesn't descend to a section of the multiplicative line bundle on $\Pi_1(T^*M)\,.$ Equivalently, the constant section is not covariantly constant along the fibers of the map 
\vspace{0.05cm}\\
\begin{equation*}
    \Pi_{\infty}^{(1)}(T^*M)\to \Pi_1^{(1)}(T^*M)\,,
    \end{equation*}
    \vspace{0.05cm}\\
with respect to the Liouville connection. To remedy this, we must define the observables using a gauge invariant, nowhere vanishing section. The existence of such a section actually implies a stronger condition than \ref{condition},  it implies that for all morphisms $X:TS^2\to T^*M,$
\vspace{0.05cm}\\
\begin{equation}\label{conditionn}
    \int_{S^2}X^*\Pi=0\,.
\end{equation}
\vspace{0.05cm}\\
A section of the line bundle $h:\Pi_{\infty}^{(1)}(T^*M)\to\mathbb{C}$ descends if and only if, for any $X:TD\to T^*M$ deforming $\gamma$
 to $\gamma',$
 \vspace{0.05cm}\\
 \begin{equation}
    h(\gamma')=h(\gamma)e^{\frac{i}{\hbar}S[X]}\,.
\end{equation}
\vspace{0.05cm}\\
Given such an section, which we assume satisfies $h\vert_M=1,$ the observables we will use instead in \ref{starry} are $\gamma\mapsto h(\gamma)t^*f(\gamma),\,\gamma\mapsto h(\gamma)s^*g(\gamma)\,.$
Using these observavles, \ref{starry} is replaced with
\vspace{0.05cm}\\
\begin{align}\label{inva}
    (f\star g)(m)=\int_{X}h(\gamma_2)t^*f(\gamma_2)\,h(\gamma_1)s^*g(\gamma_1)\,h(\gamma_3)^{-1}e^{\frac{i}{\hbar}S[X]}\,dX\,.
\end{align}    
\vspace{0.05cm}\\
The appearance of $h(\gamma_3)^{-1}$ is due to the fact that product of points in the mutiplicative line bundle over $\gamma_2,\gamma_1$ is a point over $\gamma_3\,,$ and we need to identify this point with a point over $m$ (since the integral essentially involves adding together points in different vector spaces). Now analyzing how
\vspace{0.05cm}\\
\begin{equation*}h(\gamma_2)t^*f(\gamma_2)\,h(\gamma_1)s^*g(\gamma_1)\,h(\gamma_3)^{-1}e^{\frac{i}{\hbar}S[X]}
\end{equation*}
\vspace{0.05cm}\\
changes under homotopies of the boundary of $X,$ with fixed corners, we see that the contribution from $h(\gamma_2)h(\gamma_1)h(\gamma_3)^{-1}$ exactly cancels out the contribution from $e^{\frac{i}{\hbar}S[X]}\,,$ therefore this integrand is gauge invariant. Equivalently, under the quotient map $q:\Pi_{\infty}^{(2)}(T^*M)\to \Pi_{1}^{(2)}(T^*M),$ \ref{inva} can be written as 
\vspace{0.05cm}\\
\begin{equation}
    (f\star g)(m)=\int_{X}t^*f(\gamma_2)s^*g(\gamma_1)\,e^{\frac{i}{\hbar}q^*\tilde{S}[X]}\,dX\,,
\end{equation}
\vspace{0.05cm}\\
where $\tilde{S}$ is a 2-cocycle on $\Pi_1(T^*M),$ which we will assume (as we can always do)\footnote{To understand how the section must be chosen in order to accomplish this, see the remark on page 184 of \cite{weinstein1}.}
is antisymmetric and maps to $\Pi$ under the van Est map. Using the gauge invariance, we then formally arrive at \ref{guess1}, ie.
\vspace{0.05cm}\\
\begin{equation*}
    (f\star g)(m)=\int_{s^{-1}(m)\times t^{-1}(m)} d\mu_m\otimes d\mu_m \,(t^*f)(s^*g) e^{\frac{i}{\hbar}\tilde{S}}\,.
    \end{equation*}
    \vspace{0.05cm}\\
\begin{exmp}
In the example of the Moyal product there is a canonical choice of the section $h\,.$ If $\gamma$ is a geodesic, then we define $h(\gamma)=1\,.$ This determines $h$ on all $\gamma\,.$
\end{exmp}
\section{Examples}\label{exa}$\,$
\vspace{0.05cm}\\
In this section we are going to go over some examples of \cref{guess1}. We will begin with some examples which we know deform the algebra $C_c^{\infty}(M)\,.$
\subsection{Strict Deformation Quantizations}$\,$
\begin{exmp}
Let $(M,0)$ be a manifold with the zero Poisson structure. We have already shown in the previous section that $f\star g=fg\,.$ The $C^*$-algebra is just $L^{\infty}(M,\mathbb{C})\,.$
\end{exmp}$\,$
\begin{exmp}\label{symp}
Consider the constant symplectic structure $(\mathbb{R}^2, dp\wedge dq).$ We fleshed out this example in \ref{fund}; here we will show associativity using condition \ref{ass} (this proof is adapted from \cite{Zachos}). The product is given by
\\
\begin{align}
(f\star g)(p,q)=\frac{1}{(4\pi\hbar)^2}\int_{\mathbb{R}^4}f(p'',q'')g(p',q')e^{\frac{i}{2\hbar}[(p''-p)(q-q')-(q''-q)(p-p')]}\,dp''\,dq''\,dp'\,dq'\,. 
\end{align}\\
Now by condition \ref{ass}, associativity is equivalent to the following equality, for all $(p''',q''',p'',q'',p',q',x,y)\in \mathbb{R}^6$:\\
\begin{align}
    &\int_{\mathbb{R}^2}e^{\frac{i}{2\hbar}[(p'''-x')(y'-q'')-(q'''-y')(x'-p'')+(x'-x)(y-q')-(y'-y)(x-p')]}dx'dy'\nonumber
    \\&=\int_{\mathbb{R}^2}e^{\frac{i}{2\hbar}[(p'''-x)(y-y')-(q'''-y)(x-x')+(p''-x')(y'-q')-(q''-y')(x'-p')]}dx'dy'\,.
    \end{align}\\
To see that this equality holds, we integrate over $x',y'.$ The above integrals are then equal to, respectively, \\
  \begin{align}
    &\delta(q''-q'''+y-q')\delta(p'''-p''-x+p')e^{\frac{i}{2\hbar}[-p'''q''+q'''p''
    +xq'-yp']}\,,
    \\&\delta(q'''-y+q'-q'')\delta(x-p'''+p''-p')e^{\frac{i}{2\hbar}[p'''y-q'''x-p''q'+q''p']}\,,
    \end{align}  \\
    and indeed, both of these are equal. Furthermore, the property $\overline{f\star g}=\bar{g}\star\bar{f}$ follows from antisymmetry of the cocyle and a change of variables, ie.
    \begin{align}
     \overline{f\star g}=&\frac{1}{(2\pi\hbar)^2}\int_{\mathbb{R}^4}\bar{f}(p'',q'')\bar{g}(p',q')e^{-\frac{i}{2\hbar}[(p''-p)(q-q')-(q''-q)(p-p')]}\,dp''\,dq''\,dp'\,dq'  
     \\&=\frac{1}{(2\pi\hbar)^2}\int_{\mathbb{R}^4}\bar{g}(p',q')\bar{f}(p'',q'')e^{\frac{i}{2\hbar}[(p'-p)(q-q'')-(q'-q)(p-p'')]}\,dp''\,dq''\,dp'\,dq'\\
     &=\bar{g}\star\bar{f}.
    \end{align}\\\\
We get a star product by computing the aymptotic expansion of $f\star g$ in $\hbar$ at $\hbar=0.$ To see this, use the stationary phase approximation together with integration by parts (we will flesh this out more in example \ref{good}).
\end{exmp}$\,$\\
Let us recall from the example in \ref{fund} that the Moyal product can be viewed as a gauge fixing of the associated Poisson sigma model. Related to this is the fact that on any Riemannian manifold with nonpositive sectional curvature, there is a unique geodesic in any homotopy class of paths between two points (see \cite{gallot}). This can be used to obtain the cocycle in the next two examples as well.
\vspace{0.05cm}\\
\begin{exmp}
Let's consider the symplectic 2-torus $(T^2,\omega)\,.$ In this case, using the periodicity of $f\,,g\,,$ \cref{guess1} reduces to
\vspace{0.05cm}\\
\begin{equation*}
(f\star g)(\theta,\phi)=\lim\limits_{n\to\infty}\frac{1}{(2\pi\hbar)^2}\int_{[-2\pi n,2\pi n]^4}d\theta_1d\phi_1d\theta_2d\phi_2\,f(\theta_1,\phi_1)g(\theta_2,\phi_2)e^{\frac{i}{\hbar}\omega((\theta-\theta_2,\phi-\phi_2),(\theta-\theta_1,\phi-\phi_1))}\,,
\end{equation*}
\vspace{0.05cm}\\
which is essentially the same thing as pulling back $f\,,g$ to the universal cover and using the Moyal product there. Here we are considering $(\theta,\phi)$ to be in $[0,2\pi]^2\subset\mathbb{R}^2\,.$ Indeed, one can check that 
\vspace{0.05cm}\\
\begin{equation*}
   \boxed{e^{i\theta}\star e^{i\phi}=e^{i\hbar} e^{i\theta} e^{i\phi}\,,\, e^{i\phi}\star e^{i\theta}=e^{-i\hbar} e^{i\theta} e^{i\phi}}\,,
\end{equation*}
\vspace{0.05cm}\\
so that  $e^{i\theta}\star e^{i\phi}=e^{2i\hbar}e^{i\phi}\star e^{i\theta}\,.$ Equipped with the natural $C^*$-algebra obtained by using the operator norm, this algebra is the noncommutative torus (with the involution being complex conjugation).
\end{exmp}$\,$
\vspace{0.05cm}\\
\begin{exmp}
Let $G$ be the Heisenberg group and consider the natural Poisson structure on $\mathfrak{g}^*\,.$ The symplectic groupoid is $G\ltimes\mathfrak{g}^*\rightrightarrows \mathfrak{g}^*\,.$ Note that, the exponential map $\mathfrak{g}\to G$ is a diffeomorphism and therefore we will identify $G$ and $\mathfrak{g}$ as spaces, where the latter is parameterized by $(a,b,c)$ in the usual way, ie. 
\vspace{0.05cm}\\
\begin{equation}
    \begin{pmatrix}
0 & a & c\\
0 & 0 & b\\
0 & 0 & 0
\end{pmatrix}
\end{equation}
The product can be defined via the Baker-Campbell-Hausdorff formula, ie.
\begin{equation}
    X\cdot Y=X+Y+\frac{1}{2}[X,Y]\,.
\end{equation}
On the dual basis $(x,y,z)\in\mathfrak{g}^*\,,$  the coadjoint action is given by
\begin{equation}
    (a,b,c)\cdot (x,y,z)=(x+zb,y-za,z)\,.
\end{equation}
\vspace{0.05cm}\\
Therefore, the symplectic leaves are two dimensional for $z\ne 0$ and are points for $z=0\,.$ Note that, in these coordinates the Poisson structure is given by
\vspace{0.05cm}\\
 \begin{equation}
    \Pi=\frac{1}{2}z\partial_x\wedge\partial_y\,.
\end{equation}
\vspace{0.05cm}\\
The isotropy for $z\ne 0$ is given by $\{(0,0,c):c\in\mathbb{R}\},\,$ and for $z=0$ is given by $\{(a,b,c):a,b,c\in\mathbb{R}\}\,.$
\vspace{0.5cm}\\Now that we have described the groupoid, let's describe the cocycle. There is a natural map $h:G\times \mathfrak{g}^*\to \mathbb{R}$ given
\vspace{0.05cm}\\
\begin{equation}
    (X,f)\mapsto f(X)\,,
\end{equation}
and the cocycle we are interested in is\footnote{We will explain this derivation later on.} \begin{equation}
    \delta^*h(X,Y,f)=\frac{1}{2}f[X,Y]\,.
\end{equation}
In coordinates, this gives
\begin{equation}
    \delta^*h(a',b',c',a,b,c,x,y,z)=\frac{z(ab'-ba')}{2}\,.
\end{equation}
\vspace{0.05cm}\\
By inspection, the cocycle is invariant under the whole isotropy group, therefore we can reduce the phase space by either method one or two and they give the same result. Using method 2 is a bit simpler since we quotient out by the entire isotropy, reducing $f\star g$ to an integral only involving symplectic leaves, with the canonical symplectic form. For $z\ne 0\,,$ this gives
\vspace{0.05cm}\\
\begin{align}\label{hei}
    &(f\star g)(x,y,z)=\nonumber
    \\&\boxed{\frac{1}{(2\pi z\hbar)^2}\int_{\mathbb{R}^2\times\mathbb{R}^2}f(x'',y'',z)g(x',y',z)e^{\frac{i}{z\hbar}((x''-x)(y'-y)-(x'-x)(y''-y))}\,dx''\,dy''\,dx'\,dy'}\,,
\end{align}
\vspace{0.05cm}\\
and for $z=0$ the product is $(f\star g)(x,y,0)=(fg)(x,y,0)\,.$
\vspace{0.5cm}\\This is an associative product and we obtain a $C^*$-algebra by using the operator norm, with the involution being complex conjugation. 
\vspace{0.5cm}\\Since the Heisenberg group is nilpotent, once can deformation quantize this Poisson structure according to Rieffel's prescription in \cite{rieffel2}, and though not immediately obvious, it gives the same result.
\vspace{0.5cm}\\Let us remark that, for fixed $z\,,$ \ref{hei} is the Moyal product with $\hbar\to z\hbar$ (hence in particular, it is associative). Therefore, for each symplectic leaf $\mathcal{L}_z$ we get a $\star$-representations of this algebra, which we will denote by $\mathcal{H}_{\mathcal{L}_z}\,.$ Taking the direct sum of Hilbert spaces over all symplectic leaves, ie.
\vspace{0.05cm}\\
\begin{equation}
\bigoplus_{z\in\mathbb{R}}\mathcal{H}_{\mathcal{L}_z}\,,
\end{equation}
\vspace{0.05cm}\\
we get a representation of this $C^*$-algebra. The fact that we get this representation over the direct sum of Hilbert spaces (one for each symplectic leaf) is a consequence of the $\star$-product being tangential (ie. the product restricts to symplectc leaves). Of course, classically a particle on a Poisson manifold is confined to a symplectic leaf. Therefore, a tangential product does not exhibit quantum tunneling in this sense.
\end{exmp}$\,$
\vspace{0.05cm}
\begin{exmp}
Let's generalize the example of the constant symplectic structure to get the full Moyal product. Consider a constant Poisson structure on a vector space $(V,\Pi)\,.$  The source simply connected integration is $G=V^*\times V$ with source, target and multiplication given by
\begin{equation}
\vspace{0.05cm}\\
    s(v,x)=x\,,\, t(v,x)=x+\Pi(v,\cdot)\,,\, m(v',v,x)=(v'+v,x)\,,
\end{equation}
\vspace{0.05cm}\\
and the cocycle is given by 
\vspace{0.05cm}\\
\begin{equation}
 S(v',v,x)=\pi(v,v')\,.
\end{equation}
\vspace{0.05cm}\\
Now this cocycle is invariant under the isotropy, ie. $S(v'+b,v+a,x)=S(v',v,x)\,,$ for $a\,,b$ such that $\pi(a,\cdot)=\pi(b,\cdot)=0\,.$ Therefore, we will use method 2 and quotient out by the isotropy.
The product is\footnote{Sometimes we use $2\hbar$ instead of $\hbar$ as the parameter, but it's just a convention.}
\vspace{0.05cm}\\
\begin{equation}
    (f\star g)(m)=\frac{1}{(2\pi\hbar)^{2n}}\int_{L_m\times L_m}\omega_x^n\otimes \omega_y^n\,f(x)g(y)\,e^{\frac{i}{\hbar}\Pi(x-m,y-m)}\,,
\end{equation}
\vspace{0.5cm}\\Once again, we obtain a $C^*$-algebra as discussed previously. This result agrees with the one in \cite{eli}.
\vspace{0.05cm}\\
\begin{remark}
We see that again, along with the example of the dual of the Heisenberg Lie algebra, given any morphism $\gamma\in V^*\times V$ between two objects there is a canonical morphism $\tilde{\gamma}\in \Pi_{\infty}(T^*V)\,,$ up to a path contained in the isotropy,  and our cocycle $S$ is the pullback of the canonical one.
\end{remark}
\end{exmp}$\,$
\vspace{0.05cm}
\begin{exmp}\label{good}
Let's consider one final example. We will once again consider $(\mathbb{R}^2,dx\wedge dy),$ but we will use a different cocycle. In fact, in this example we will use one that isn't antisymmetric. We will see that we still get a deformation quantization, however it will have a property that the previous examples didn't. \\\\
This time we will describe the groupoid by $\text{Pair}(\mathbb{R}^2)\rightrightarrows \mathbb{R}^2.$ There is a unique arrow between any two objects, therefore we will specify an arrow by $(x',y',x,y)$ and a composable pair by $(x'',y'',x,y,x',y').$
Consider 
\begin{equation}
    h(x',y',x,y)=\varepsilon y'y,
    \end{equation}
for $\varepsilon\in\mathbb{R}.$ We have that $d_{\text{CE}}VE(h)=0,$ and since $VE$ is a cochain map (see \ref{appen}) it follows that $VE(\delta^*h)=0.$ Explicitly,
\begin{equation}
 \delta^*h(x'',y'',x,y,x',y')=\varepsilon(yy'-y'y''+yy'').   
\end{equation}
Now consider the cocycle\\
\begin{equation}
    S[x'',y'',x,y,x',y']=(x''-x)(y-y')-(y''-y)(x-x')+\varepsilon(yy'-y'y''+yy'').
\end{equation}\\
This is the same as the cocycle in \ref{fund}, with the addition of $\delta^*h,$ therefore $VE(S)=\partial_x\wedge\partial_y.$ \\\\
Now the first step is to normalize, ie. we want $1\star 1=1.$ We have that \\
\begin{equation}
    1\star 1=\int_{\mathbb{R}^2\times\mathbb{R}^2}e^{\frac{i}{\hbar}\big[(x''-x)(y-y')-(y''-y)(x-x')+\varepsilon(yy'-y'y''+yy'')\big]}\,dx''dy''dx'dy'\,.
\end{equation}
First doing the integrals over $x'',x'$ we get\\
\begin{align}
&(2\pi\hbar)^2\int_{\mathbb{R}\times\mathbb{R}}\delta(y-y')\delta(y''-y)e^{\frac{i}{\hbar}\big[-x(y-y')-(y''-y)x+\varepsilon(yy'-y'y''+yy'')\big]}\,dy''dy'\\
&=(2\pi\hbar)^2e^{\frac{i}{\hbar}\varepsilon y^2}\,.
\end{align}\\
Therefore, \ref{guess1} is
\begin{align}
    &(f\star g)(x,y)=
    \\&\frac{1}{(2\pi\hbar)^2}\int_{\mathbb{R}^2\times\mathbb{R}^2}f(x'',y'')g(x',y')e^{\frac{i}{\hbar}\big[(x''-x)(y-y')-(y''-y)(x-x')+\varepsilon(yy'-y'y''+yy''-y^2)\big]}\,dx''dy''dx'dy'\,.\nonumber
\end{align}\\
Also, by the stationary phase approximation (see \cite{bates}) \\
\begin{equation}
    \lim\limits_{\hbar\to 0}(f\star g)(m)=f(m)g(m).
\end{equation}$\,$
\\\\Now we will show that the this product gives a deformation quantization of $\displaystyle\partial_x\wedge \partial_y.$ Consider $\frac{d}{d\hbar}(f\star g)(m).$ After integrating by parts, we find that it is equal to\\
\begin{align}\label{equu}
    &\frac{-i}{(2\pi\hbar)^2}\int_{\mathbb{R}^2\times\mathbb{R}^2}\big(\partial_{x''}\partial_{y'}-\partial_{x'}\partial_{y''}\big)f(x'',y'')g(x',y')e^{\frac{i}{\hbar}S[x'',y'',x,y,x',y']}\,dx''dy''dx'dy'\nonumber
    \\&-\frac{2i\varepsilon}{(2\pi)^2\hbar^4}\int_{\mathbb{R}^2\times\mathbb{R}^2}f(x'',y'')g(x',y')\,(y''-y)(y'-y)e^{\frac{i}{\hbar}S[x'',y'',x,y,x',y']}\,dx''dy''dx'dy'\,.
\end{align}\\
Now by the stationary phase approximation, the first integral in \ref{equu} converges to $i\{f,g\}$ as $\hbar\to 0,$ so we only have to compute the second one. Let's make the following change of variables:\\
\begin{align}
Y''=y''-y,\;Y'=y'-y,\;X''=x''-x,\; X'=x'-x.
\end{align}\\
Then our integral becomes\footnote{We are supressing the arguments of $f,g$ for brevity.}\\
\begin{align}
    -\frac{2i\varepsilon}{(2\pi)^2\hbar^4}\int_{\mathbb{R}^2\times\mathbb{R}^2}fg\,Y'Y''e^{\frac{i}{\hbar}[-X''Y'+Y''X'-Y'Y'']}\,dX''dY''dX'dY'\,.
\end{align}\\\\
Now the exponent $-X''Y'+Y''X'-Y'Y''$ has a nondegenerate and unique critical point at $(X'',Y'',X',Y')=(0,0,0,0),$ so the idea is to use the stationary phase approximation. A priori the stationary phase approximation only tells us that the integral $\sim (2\pi\hbar)^2$ as $\hbar\to 0,$ which isn't enough since we have an $\hbar^4$ in the denominator. Inspired by the Morse lemma we make the following substitution:\footnote{And thank you to Willie Wong for help with this \cite{willie}.}\\
\begin{align}
    A=\frac{X''+Y'}{2},\; B=\frac{X''-Y'}{2},\; C=\frac{Y''+X'-y}{2},\; D=\frac{Y''-X'+y}{2}\,.
\end{align}\\
Then our integral becomes \\
\begin{align}
     \frac{2i\varepsilon}{(2\pi)^2\hbar^4}\int_{\mathbb{R}^2\times\mathbb{R}^2}fg\,(A-B)(C+D)e^{\frac{i}{\hbar}[-A^2+B^2+C^2-D^2]}\,dA\,dB\,dC\,dD\,.
\end{align}\\
Now, after integrating by parts, we find that this is equal to\\
\begin{align}
     -\frac{i\varepsilon}{2(2\pi\hbar)^2}\int_{\mathbb{R}^2\times\mathbb{R}^2}[(\partial_A+\partial_B)(\partial_C-\partial_D)fg]\,e^{\frac{i}{\hbar}[-A^2+B^2+C^2-D^3]}\,dA\,dB\,dC\,dD\,,
\end{align}\\
and as $\hbar\to 0$ this goes to $-i\frac{\varepsilon}{2}\,\partial_yf(x,y)\partial_yg(x,y).$
Therefore, we have that \\\\
\begin{equation}
    (f\star g)(x,y)\sim f(x,y)g(x,y)-i\hbar\big[\{f,g\}(x,y)+\frac{\varepsilon}{2}\,\partial_yf(x,y)\partial_y g(x,y)\big]\;\;\;(\hbar\to 0)\,.
\end{equation}$\,$\\
From which it follows that \\
\begin{equation}\label{true}
    f\star g-g\star f\sim -2i\hbar\{f,g\}\;\;\;(\hbar\to 0)\,.
\end{equation}\\
Now we see the difference between this example and the previous examples: the first order term of $f\star g$ in $\hbar$ isn't just the Poisson bracket, it contains a symmetric part. As indicate by Kontsevich in \cite{kontsevich} (page 3), the symmetric part can always be killed by a gauge transformation. In this case the gauge transformation is 
\begin{equation}f\to f-i\hbar\frac{\epsilon}{4}\partial_y^2,
\end{equation}
since $-f\partial_y^2g+\partial_y^2(fg)-g\partial_y^2 f=2\partial_yf\,\partial_y g.$
\\\\Indeed, the symmetric part of \ref{true} corresponds to the symmetric part of our cocycle, $\varepsilon(yy'-y'y''+yy'')\,.$ See \ref{cocycle} for a brief discussion on this.
\end{exmp}$\,$
\vspace{0.05cm}\\
\subsection{Lie-Poisson Structures of Exponential and Semisimple Lie Groups}$\,$
\vspace{0.05cm}\\
For the Poisson structures associated to exponential Lie groups (ie. Lie groups for which $\text{exp}:\mathfrak{g}\to G$ is a diffeomorphism) and semisimple Lie groups, there is canonical degree 2-cocycle on the symplectic groupoid $G\ltimes \mathfrak{g}^*\rightrightarrows \mathfrak{g}^*$ with the desired properties, which we will now describe. First note that, there is a canonical Lie algebroid one form on $T^*\mathfrak{g}^*\cong \mathfrak{g}^*\times \mathfrak{g}\,,$ given by 
\begin{equation}
    \alpha(f,v)=f(v)\,.
\end{equation}
The significance of this one form is that $d_{CE}\alpha$ is the Poisson structure, where $d_{CE}$ is the Chevalley-Eilenberg differential.
\vspace{0.05cm}\\
\begin{exmp}\label{expp}
Let $G$ be an exponential Lie group. Then there is a canonical function on $G\times \mathfrak{g}^*\,,$ given by
\begin{equation}
    \tilde{\alpha}(g,f)=f(\log{g})\,.
\end{equation}
Now we have that $VE(\tilde{\alpha})=\alpha\,,$ and since $VE$ is a cochain map, it follows that $VE(\delta^*\tilde{\alpha})$ is the Poisson structure. It follows that the antisymmetrization of $\delta^*\tilde{\alpha}$ has the desired properties and we can take this to be the cocycle.
\vspace{0.5cm}\\A particularly nice class of examples of exponential groups are nilpotent Lie groups. These are nice because the Baker-Campbell-Hausdorff formula has only finitely many terms and we can identify the exponential map with the identity map on the Lie algebra, in which case the Haar measure is just the Lebesgue measure on the underlying vector space.
\end{exmp}$\,$
\vspace{0.05cm}\\
\begin{exmp}\label{sphere}
Let $G$ be a semisimple Lie group. We can use the Killing form, given by
\begin{equation}
    \mathfrak{g}\times\mathfrak{g}\to\mathbb{R}\,, (X,Y)\mapsto \text{tr}(\text{Ad}_X\text{Ad}_Y)\,,
\end{equation}to identify $\mathfrak{g}^*$ with $\mathfrak{g}\,,$ hence we identify the Lie algebroid with $\mathfrak{g}\times\mathfrak{g}$ and we can identify the symplectic groupoid with $G\ltimes \mathfrak{g}\rightrightarrows \mathfrak{g}\,.$ The canonical one form $\alpha$ is then described as 
\begin{equation}
    \alpha_X(Y)=\text{tr}(\text{Ad}_X\text{Ad}_Y)\,.
\end{equation}
On $G\times\mathfrak{g}$ there is a natural function given by \begin{equation}\label{semi}
    \tilde{\alpha}(g,X)=\text{tr}(\text{Ad}_g\text{Ad}_X)\,,
\end{equation}
and $VE(\tilde{\alpha})=\alpha\,,$ so as in the previous example the antisymmetrization of $\delta^*\tilde{\alpha}$ is our desired cocycle.
\vspace{0.5cm}\\The nicest case of this example is when $G$ is compact, since then there are no divergence issues. We can write $\star$ in the following way\footnote{The only thing that needs to be checked is that $C_{X,\hbar}\ne 0$ for any $X\,.$ At the very least we know that $C_{X,\hbar}\ne 0$ near $X=0\,.$}:
\begin{equation}\label{compact semisimple}
\boxed{(f\star k)(X)=C_{X,\hbar}\int_{G\times G} dg\,dg' \,f(\text{Ad}_g \,X)k(\text{Ad}_{g'}\, X)e^{\frac{i}{\hbar}S(g,g',X)}}\,,
\end{equation}
where
\begin{align}
   &S(g,g',X)=\frac{1}{2}\,[\delta^*\tilde{\alpha}(g,g',X)-\delta^*\tilde{\alpha}(g'^{-1},g^{-1},\text{Ad}_{gg'}X)]\,,\\
    & \tilde{\alpha}(g,X)=\text{tr}(\text{Ad}_g\text{Ad}_X)\,.
    \end{align}
Here $dg$ is the normalized Haar measure on $G$ and $C_{X,\hbar}$ is the normalization parameter, which in this case is constant along symplectic leaves (ie. it is a Casimir function), and $C_{0,\hbar}=1\,.$
\vspace{0.5cm}\\This product satisfies
\begin{equation}
    \overline{(f\star k)}=\bar{k}\star\bar{f}
\end{equation}
and this product is $G$-equivariant in the sense that 
\begin{equation}
    (g\cdot f)\star(g\cdot k)=g\cdot(f\star k)\,.
\end{equation}
Therefore, the $C^*$-algebra associated with $\star$ comes with a $G$-action, since
\begin{equation}
    (g\cdot f)\star k=g\cdot(f\star (g^{-1}\cdot k))\,.
\end{equation}
This $G$-action respects the involution and the norm.
\vspace{0.5cm}\\Specializing to $G=SU(2)$ and restricting $\star$ to a generic symplectic leaf, we get an $SU(2)$-equivariant map from $L^2(S^2)$ into the restricted $C^*$-algebra. Results by Rieffel (see the next subsection) suggest that $\star$ is nonassociative, so this (probably) at best constitutes a strict quantization of $C^{\infty}(S^2)$ (but not a strict deformation quantization).\\
\vspace{0.05cm}\\Note that, for semisimple subgroups of $GL_n(\mathbb{C})$ there are simpler options for the cocycle which can be obtained by rewriting the Killing form. For example the Killing form on $\mathfrak{s}\mathfrak{u}(n)$ is equal to $2n\text{tr}(XY),$ therefore instead of \ref{semi} one can use $\tilde{\alpha}(g,X)=\text{tr}(gX)\,.$
\end{exmp}
\subsubsection{Associativity}
Before concluding we should discuss associativity of examples \ref{expp} and \ref{sphere}. We have seen one example where \ref{exp} is associative (the Heisenberg group example). However, in general it seems doubtful that this is the case since in \cite{gamme} it is shown that some nilpotent Lie groups do not have tangential star products. However, there is a possible execption: if the cocycle has a critical point away from the identity bisection, then the asymptotic expansions in $\hbar$ at $\hbar=0$ need not give a star product in the sense of \cite{kontsevich}. \\\\
In the case of the duals semisimple Lie algebras, a similar obstruction to the existence of tangential star products was shown in \cite{cahen}. Furthermore, in \cite{rieffel3} (in example 14) it is stated that there is no strict deformation quantization of symplectic $S^2$ which preserves the action of $SO(3).$ Of course, this doesn't mean that there isn't a \textit{strict quantization} (see \ref{explain}) of symplectic $S^2$ which preserves the action of $SO(3).$\\\\
In the case of regular Poisson manifolds, the existence of a formal tangential product (but not necessarily a star product) is proved in \cite{Masmoudi} (note that the dual of the Heisenberg Lie algebra isn't regular, so it's not an if and only if). Related to all of this is the conjecture of Weinstein \ref{weinc}.
    \section{Final Remarks}
To recap, based on work done by Bonechi, Cattaneo, Felder and Zabzine on Poisson sigma models, we formally showed that Kontsevich's star product can be obtained from the geometric quantization of a Lie 2-groupoid. We then interpreted the Poisson sigma model as a sigma model valued in higher groupoids and argued that, if one can perform symplectic reduction (including dealing with infinities) on the quantization, then one can use a (non-perturbative) groupoid-valued sigma model to compute correlation functions of gauge invariant observables.\footnote{Again, our notion of gauge equivalence is coarser than the one in \cite{bon}.}\\\\
One thing that isn't clear is how to precisely connect this work with geometric quantization, in the sense of \cite{eli}. The geometric quantization of the Lie-2 groupoid, as well as the groupoid-valued sigma model, each have many things in common with the usual geometric quantizaion of Poisson manifolds, but neither are the same — it's not completely clear how polarizations enter the picture, though they are used for reduction, in a sense. \\\\
Now even in cases where the products aren't associative we still get $C^*$-algebras, which could conceivably form strict quantizations.\footnote{See \ref{explain}.} Therefore, even in cases where a tangential deformation quantization doesn't exist, it is still conceivable that a tangential strict quantization does. Some possible examples come from compact semisimple Lie groups.\\\\
Let's examine all of this from the perspective of the Poisson sigma model and the original product given by \ref{star}. Consider the example of $S^2$ with its $SO(3)$-invariant symplectic structure. Rieffel proved that there is no strict deformation quantization of $S^2$ which preserves the action of $SO(3)$ (see \cite{rieffel3}). If a \textit{non-perturbative} definition of the Poisson sigma model exists in this case then the author would expect that it would induce a strict deformation quantization of $S^2$ by using the operator norm. However, since it couldn't be $SO(3)$-invariant, this could mean that either \ref{star} isn't actually $SO(3)$-invariant (despite the Poisson structure being), or the product is not associative, or the product is not continuous with respect to the operator norm. Note that, Kontsevich's proof of associativity of the star product in \cite{kontsevich} only shows that \ref{star} is \textit{perturbatively associative} when expanding about the constant critical point. \\\\
Now let's discuss some generalizations of the ideas in this paper. First of all, it isn't necessary to use the source simply connected integration in the groupoid-valued sigma model, we only used it because if the 2-cocycle exists on any 1-groupoid then it will exist on the source simply connected one. However for the duals of Lie algebras, for example, we can use any integration. In addition, requiring the cocycle to be antisymmetric isn't necessary either. Furthermore, we didn't use the smooth structure on the groupoid much (mostly just in some neighborhood of the identity bisection), so there is a possibility that something similar can be carried out even in the case where the Poisson structure isn't integrable.\\\\
Lastly, in the case where the condition \ref{conditionn} is violated, there are two options one could try. One is to embed the Poisson manifold as a closed submanifold inside another and try to obtain a quantization by quantizing the ambient Poisson manifold, as we did in example \ref{sphere}. Another option is to make all of the constructions in this paper using the integration of $T^*M\oplus_{\Pi} \mathbb{R}\to M,$ ie. rather than applying these constructions to a groupoid integrating $T^*M\to M,$ we apply them to a groupoid integrating $T^*M\oplus_{\Pi} \mathbb{R}\to M,$ the natural extension of $T^*M$:\\
\begin{equation}
0_M\to \mathbb{R}\to T^*M\oplus_{\Pi} \mathbb{R}\xrightarrow[]{p} T^*M\to 0_M\,.
\end{equation}$\,$\\
The only difference is that, rather than looking for a cocycle $S$ on $\Pi_1(T^*M)$ such that $VE(S)=\Pi,$  we look for a cocycle $\tilde{S}$ on $\Pi_1(T^*M\oplus_{\Pi} \mathbb{R})$ such that $VE(\tilde{S})=p^*\Pi.$ Such a cocycle will always exist (as long as $\Pi_1(T^*M\oplus_{\Pi} \mathbb{R})$ exists) since $p^*\Pi$ is (canonically) trivial in cohomology. If the cocycle satisfies the required condition used in \cite{bon}, then $T^*M\oplus_{\Pi} \mathbb{R}$ is integrable.
\subsection*{Acknowledgements}
During the writing of this paper I exchanged several helpful emails with Alberto S. Cattaneo.
    \section{Appendix}\label{appen}
Here we are going to briefly discuss what we need to know about cohomology and the van Est map for this paper. Everything we do here works in broader contexts but we will state just what we need for this paper. For references about the topics here, see \cite{Crainic}, \cite{eli}, \cite{Lackman}, \cite{Lackman2}, \cite{Meinrenken}, \cite{weinstein1}.
\subsection{Cohomology of Simplicial Spaces}$\,$
\vspace{0.05cm}\\
Let $S^{\bullet}$ be a simplicial space. There is a differential $\delta^*:C(S^{(n)},\mathbb{R})\to C(S^{(n+1)},\mathbb{R})$ given by 
\begin{equation}
    \delta^*f=\sum_{i=0}^n (-1)^id_i^*f\,,
\end{equation}
where $d_i:S^{(n+1)}\to S^{(n)},\,i=0,\ldots,n,$ are the corresponding face maps. Being a differential, $\delta^*$ satisfies $\delta^{*2}=0\,,$ and therefore we get a cochain complex
\begin{equation}
   C(S^{(0)},\mathbb{R})\to C(S^{(1)},\mathbb{R})\to C(S^{(2)},\mathbb{R})\to\cdots\,.
\end{equation}
\begin{definition}(see \cite{Deligne})
The cohomology of this simplicial space, with coefficients in the sheaf $\mathcal{O}\,,$ is defined to be the cohomology of this cochain complex, and is denoted $H^*(S^{\bullet},\mathbb{R})\,.$ When $S^{\bullet}$ is the nerve of a groupoid $G\,,$ this is called groupoid cohomology with coefficients in $\mathcal{O}\,;$ we will denote this by $H^*(G^{\bullet},\mathbb{R})\,.$
\end{definition}
Note that, if the simplicial space is smooth then it is understood that we only consider smooth functions in each degree. Furthermore, we can use $\mathbb{C}$ instead of $\mathbb{R}$ in the cochain complex without modifying any of the formulas. In fact, we can consider much more general coefficients for cohomology, like the sheaf of $S^1$-valued functions, but the definition is a bit more complicated since we need to take a resolution of the sheaf in each degree, see \cite{Deligne}, \cite{Lackman}.
\begin{proposition}(see \cite{eli})
Associated to any $2$-cocycle $f$ on a Lie groupoid $G\rightrightarrows M$ is a central extension of Lie groupoids
\vspace{0.05cm}\\
\begin{equation}
    M\to S^1_M\to S^1\times_{f}G\to G\to M\,,
\end{equation}
\vspace{0.05cm}\\
where $M$ denotes the Lie groupoid with only identity morphisms and $S^1_M$ denotes $M\times S^1\,,$ thought of as a Lie groupoid over $M$ with the source and target maps being the projection and the composition being the one induced by $S^1\,.$ Here $S^1_M\times_{f}G$ is the groupoid whose objects are $M$ and whose space of arrows is $S^1\times G^{(1)}\,,$ with the composition given by $(\lambda,g_2)\cdot(\beta,g_1)=(\lambda\beta e^{if(g_2,g_1)},g_2g_1)\,.$\\
\vspace{0.05cm}\\Furthermore, associated to such an extension is a multiplicative complex line bundle (ie. a line bundle over $G^{(1)}$ with a compatible multiplication), whose total space is $\mathbb{C}\times G^{(1)}\,,$ and with a composition given by the same formula as above, ie. $(\lambda,g_2)\cdot(\beta,g_1)=(\lambda\beta\, e^{if(g_2,g_1)},g_2g_1)\,.$
More generally, we have the following:
\vspace{0.05cm}\\
\begin{definition}(see \cite{eli}, \cite{weinstein1}
A mutliplicative line bundle $L\to G$ is given by a line bundle $L\to G^{(1)}$ together with a nowhere vanishing linear map
\begin{equation}\label{induce}
    s:d_0^*L\otimes d_2^*L\to d_1^*L\,,
\end{equation}
such that, for a composable pair $(g_2,g_1),$ the induced map 
\begin{equation}
    L_{g_2}\times L_{g_1}\to  L_{g_2}\otimes L_{g_1}\xrightarrow[]{s} L_{g_2g_1}
    \end{equation}
is associative. Here $L_{g}$ means points in the line bundle over $g\,.$
\end{definition}
\end{proposition}
\subsection{Cohomology of Lie Algebroids}
\begin{definition}(\cite{Lackman}, \cite{Lackman2}, \cite{weinstein1})
Let $\pi:\mathfrak{g}\to M$ be a Lie algebroid with anchor map $\alpha\,,$ and let $\mathcal{C}^n(\mathfrak{g},\mathbb{R})$ denote the real-valued $n$-forms on $\mathfrak{g}\,.$ There is a canonical differential given by
 \begin{align*}
     d_\text{CE}:\mathcal{C}^n(\mathfrak{g},\mathbb{R})\to \mathcal{C}^{n+1}(\mathfrak{g},\mathbb{R})\,,\,n\ge0
 \end{align*}
 defined as follows: let $\omega\in \mathcal{C}^n(\mathfrak{g},\mathbb{R})$ and pick an $m\in M;$ we want to compute $d_{CE}\omega$ at the point $m\,.$ Then for $X_1\,,\ldots \,, X_{n+1}\in \pi^{-1}(m)\,,$ choose local extensions $\mathbf{X}_1\,,\ldots\,,\mathbf{X}_{n+1}$ of these vectors, ie. 
 \begin{align*}
     p\mapsto\mathbf{X}_1(p)\,,\ldots\,,p\mapsto\mathbf{X}_{n+1}(p)
     \end{align*}
     are sections of $\mathfrak{g}$ over some open set $U$ containing $m\,,$ and such that $\mathbf{X}_i(m)=X_i$ for all $1\le i\le n+1)\,.$\\\\ Then let
 \begin{align}
     d_{\text{CE}}\omega(X_1\,,\ldots\,,X_{n+1})&=\sum_{i<j}(-1)^{i+j-1}\omega([\mathbf{X}_i,\mathbf{X}_j],\mathbf{X}_1,\ldots ,\hat{\mathbf{X}}_i,\ldots
     , \hat{\mathbf{X}}_j,\ldots , \mathbf{X}_{n+1})
\vert_{p=m}     \\& +\sum_{i=1}^{n+1}(-1)^i {\alpha(\mathbf{X}_i)}(\omega(\mathbf{X}_i,\ldots , \hat{\mathbf{X}}_i ,\ldots , \mathbf{X}_{n+1}))\vert_{p=m}\,,
 \end{align}
 This is well-defined and independent of the chosen extensions. We will denote this cohomology by $H^*(\mathfrak{g},\mathbb{R})\,.$
\end{definition}
Note that, as in the case of cohomology of simplicial spaces, we can replace $\mathbb{R}$ with $\mathbb{C}$ and all the formulas remain the same. Furthermore, we can take the cohomology of Lie algebroids with respect to more general sheaves, like the sheaf of $S^1$-valued functions, which requires a resolution to be performed, see \cite{Lackman}.
\begin{proposition}(see \cite{luk})
Associated to any closed Lie algebroid $2$-form $\omega$ is a central extension of Lie algebroids\\
\begin{equation}
    0_M\to \mathbb{R}_M\to \mathbb{R}_M\oplus_{\omega}\mathfrak{g}\to\mathfrak{g}\to 0_M\,,
\end{equation}$\,$\\
where $0_M$ denotes the zero vector bundle over $M$ and $\mathbb{R}_M$ denotes the vector bundle $M\times\mathbb{R}\,,$ thought of as a Lie algebroids. As a vector bundle, $\mathbb{R}_M\oplus_{\omega}\mathfrak{g}$ is just $\mathbb{R}_M\oplus\mathfrak{g}\,.$ The Lie algebroid $\mathbb{R}_M\oplus_{\omega}\mathfrak{g}$ has anchor map given by $(\lambda,X)\mapsto\alpha(X)\,,$ and the bracket is given by
\begin{equation}
    [(\lambda,X),(\beta,Y)]=(\alpha(X)\beta-\alpha(Y)\lambda+\omega(X,Y),[X,Y])\,.
\end{equation}$\,$
\vspace{0.05cm}\\
Two such extensions are isomorphic if and only if the $2$-forms agree in cohomology.
\end{proposition}
\subsection{The van Est map}
The simplest\footnote{The van Est map works for more general coefficients than $\mathbb{R}$ or $\mathbb{C},$ in particular you can take them in any abelian Lie group, see \cite{Lackman}. Furthermore, the van Est map also works on differentiable stacks, see \cite{Lackman2}.} manifestation of the van Est map is as a map from the cohomology of a Lie groupoid to the cohomology of its Lie algebroid, and is denoted 
\vspace{0.05cm}\\
\begin{equation*}
    VE:H^*(G^{\bullet},\mathbb{R})\to H^*(\mathfrak{g},\mathbb{R})\,.
\end{equation*}
\vspace{0.05cm}\\However, it is important to note that the van Est map exists at the level of cochains, as we use this fact in the paper. Before defining it, we need some notation. \\
\vspace{0.05cm}\\Let $S_n$ denotes permutations on $\{1,2,\ldots, n\}$ and for $\sigma\in S_n$ let $\text{sgn}(\sigma)$ be the sign of the permutation. Furthermore, given $X\in\Gamma(\mathfrak{g})\,,$ we can right translate it to a vector field $R_X$ on $G^{(1)}\,.$  Now suppose we have a function $f:G^{(n)}\to\mathbb{R}\,,$ $n\ge 1\,.$ We get a new function $R_Xf:G^{(n-1)}\to\mathbb{R}$ by defining
\vspace{0.05cm}\\
\begin{equation}
    R_Xf(g_{n-1},\ldots, g_1)=R_Xf(g_{n-1},\ldots, g_1,\cdot)\vert_{s(g_1)}\,,
\end{equation}$\,$
\vspace{0.05cm}\\
ie. we differentiate it in the first component (or last, depending on your point of view) at the identity $s(g_1)\,.$ Now we are ready to define the van Est map, at the level of cochains.
\vspace{0.05cm}\\
\begin{definition}(see \cite{Crainic}, \cite{Lackman}, \cite{Lackman2}, \cite{weinstein1})
Let $f:G^{(n)}\to \mathbb{R}$ and let $X_1,\ldots X_n\in\Gamma(\mathfrak{g})\,.$ Then we define
\vspace{0.05cm}\\
\begin{align}
    VE(f)(X_1,\ldots, X_n)=\sum_{\sigma\in S_n} \text{sgn}(\sigma)R_{X_{\sigma(n)}}\cdots R_{X_{\sigma(1)}}f\,,
\end{align} $\,$
\vspace{0.05cm}\\
This definition makes sense even if $X_1,\ldots,X_n$ are only defined over a point $m\in M,$\footnote{We only assumed they were sections of $\mathfrak{g}$ because it's slightly quicker to define, however looking at the formula one can see that it isn't even necessary to choose extensions.} and this defines our $n$-form on $\mathfrak{g}\,.$
\end{definition}$\,$
\vspace{0.05cm}\\
Now the following lemma tells us that $VE$ descends to a map on cohomology:\\
\begin{lemma}(see \cite {weinstein})
$VE$ is a cochain map. 
\end{lemma}$\,$
\vspace{0.05cm}\\
An important theorem, called the van Est isomorphism theorem, is the following (originally proved by Crainic in \cite{Crainic}):
\begin{theorem}\label{VE}
Suppose $G$ has $n$-connected source fibers. Then $VE:H^*(G,\mathbb{R})\to H^*(\mathfrak{g},\mathbb{R})$ is an isomorphism up to and including degree $n\,,$ and is injective in degree $n+1\,.$ Its image in degree $n+1$ consists of those cohomology classes $[\omega]$ such that the right translation of $\omega$ to each source fiber is exact.
\end{theorem}$\,$
\vspace{0.05cm}\\
For this paper, the only cases we will be concerned with are degrees $1$ and $2$ (but mostly $2$). The van Est isomorphism theorem will tell us when we can expect to integate the 2-cocycle $\Pi$ (associated to a Poisson structure) to a 2-cocycle on the Lie groupoid. In particular, this is the case if the source fibers are 2-connected (ie. if they have vanishing homotopy groups up to and including degee 2).


\begin{thebibliography}{9}
\bibitem{amar}
N. Ben Amar. \textit{Tangential deformations on the dual of nilpotent special Lie algebras.} Pacific
J. Math. 170 (1995), 297–318.
\bibitem{baker}
George A. Baker. \textit{Formulation of Quantum Mechanics Based on the Quasi-Probability Distribution Induced on Phase Space.} Jr. Phys. Rev. 109, 2198 – Published 15 March 1958 (doi:10.1103/PhysRev.109.2198)
\bibitem{bates}
Sean Bates, Alan Weinstein. \textit{Lectures on the geometry of quantization.} Berkeley Mathematics Lecture Notes
Volume: 8; 1997; 137 pp.
\bibitem{Pierre}
Pierre Bieliavsky and Victor Gayral. \textit{Deformation Quantization for Actions of Kählerian Lie Groups.} Memoirs of the American Mathematical Society, Volume: 236; 2015; 154 pp.
\bibitem{bojo}
M. Bojowald, A. Kotov, and T. Strobl. \textit{Lie algebroid morphisms, Poisson sigma models, and
off-shell closed gauge symmetries.} J. Geom. Phys. 54 (2005) 400–426.
 \bibitem{bon}
F. Bonechi, A. S. Cattaneo and M. Zabzine, \textit{Geometric quantization and
non-perturbative Poisson sigma model.} Adv. Theor. Math. Phys. 10 (2006)
683 [arXiv:math/0507223].
\bibitem{zab}
F. Bonechi and M. Zabzine. \textit{Lie algebroids, Lie groupoids and TFT.} Journal of Geometry and Physics, vol. 57, 731-744 (2005).
\bibitem{bu}
M. Buneci. \textit{Groupoid C*-algebras.} Surveys in Mathematics and its Applications, 1 (2006), 71-98. MR2274294(2007h:46067). Zbl 1130.22002.
\bibitem{ale}
A. Cabrera. \textit{Generating Functions for Local Symplectic Groupoids and Non-perturbative Semiclassical Quantization.} Commun. Math. Phys. 395, 1243–1296 (2022).
\bibitem{cahen}
 M. Cahen, S. Gutt and J. Rawnsley. \textit{On tangential products for the coadjoint Poisson structure.} Commun.Math. Phys. 180, 99–108 (1996). https://doi.org/10.1007/BF02101183
 \bibitem{catt}
 A. S. Cattaneo and G. Felder. \textit{A Path Integral Approach to the Kontsevich Quantization Formula}. Comm Math Phys 212, 591–611 (2000). https://doi.org/10.1007/s002200000229
 \bibitem{felder}
A.S. Cattaneo and G. Felder. \textit{Poisson sigma models and symplectic groupoids.} Quantization of singular symplectic quotients, Progr. Math., vol. 198,
Birkh¨auser, Basel, 2001, pp. 61–93.
 \bibitem{Crainic} 
Marius Crainic.
\textit{Differentiable and Algebroid Cohomology, van Est Isomorphisms, and Characteristic Classes.} 
Commentarii Mathematici Helvetici, Vol.78, (2003) pp. 681-721.
\bibitem{rui}
Marius Crainic and Rui Loja Fernandes.
\textit{Integrabiltiy of Lie Brackets.}
Annals of Mathematics, 157 (2003), 575-620.
\bibitem{ruif}
Marius Crainic and Rui Loja Fernandes. \textit{Integrability of Poisson Brackets.} J. Differential Geom. 66 (1) 71 - 137, January, 2004. https://doi.org/10.4310/jdg/1090415030
\bibitem{zhu}
Marius Crainic and Chenchang Zhu.
\textit{Integrability of Jacobi and Poisson structures.}
Annales de l'Institut Fourier, Volume 57 (2007) no. 4, pp. 1181-1216.
\bibitem{Curtright}
Thomas L. Curtright, David B Fairlie and and Cosmas K. Zachos. \textit{ Concise Treatise on Quantum Mechanics in Phase Space.} Default Book Series (2014).
Default Book Series. January 2014
\bibitem{Deligne} 
P. Deligne.
\textit{Th\'{e}orie de Hodge. III.} 
Inst. Hautes \'{E}tudes Sci. Publ. Math No. $\mathbf{44}$ (1974),
5-77.
\bibitem{gallot}
Sylvestre Gallot, Dominique Hulin, Jacques Lafontaine. \textit{Riemannian Geometry.} Springer-Verlag Berlin Heidelberg 2004.
 \bibitem{gamme}
 A. Gammella. \textit{Tangential products.} Letters in Mathematical Physics 51 , 1–15 (2000).
 \bibitem{getzler}
Ezra Getzler. \textit{Lie theory for nilpotent L$\infty$ algebras.} Annals of Mathematics, 170 (2009), 271–301
\bibitem{stern}
V. Guillemin, S. Sternberg. \textit{Geometric quantization and multiplicities of group representations.} Inventiones Mathematicae (1982), 67 (3): 515–538.
\bibitem{Gleason}
Jonathan Gleason. \textit{The C*-algebraic formulation of quantum mechanics}. (2009) \href{https://math.uchicago.edu/~may/VIGRE/VIGRE2009/REUPapers/Gleason.pdf}{THE C*-ALGEBRAIC FORMALISM OF QUANTUM MECHANICS}
\bibitem{Gorb}
V.V. Gorbatsevich, A.L. Onishchik. \textit{Lie Groups and Lie Algebras II.} Springer, Feb 3, 2000. 
\bibitem{luk}
Marco Gualtieri and Kevin Luk. \textit{Log Picard Algebroids and Meromorphic Line Bundles.} International Mathematics Research Notices, Volume 2021, Issue 21, November 2021, Pages 16592–16635.
 \bibitem{hall}
 B.C. Hall. \textit{Holomorphic methods in analysis
and mathematical physics.} Contemp. Math.
260 (2000). 
 \bibitem{eli}
Eli Hawkins. \textit{A Groupoid Approach to Quantization.}
J. Symplectic Geom. 6 (2008), no. 1, 61-125.
\bibitem{andre}
André Henriques.
\textit{Integrating L$\infty$-algebras.}
Compositio Mathematica (2008), Vol 144.
\bibitem{heinz}
M. Heins, O. Roth, S. Waldmann. \textit{Convergent products on cotangent bundles of Lie groups.} Math. Ann. (2022). https://doi.org/10.1007/s00208-022-02384-x
\bibitem{kontsevich}
M. Kontsevich. \textit{Deformation Quantization of Poisson Manifolds.} Letters in Mathematical Physics 66, 157–216 (2003). https://doi.org/10.1023/B:MATH.0000027508.00421.bf
\bibitem{Lackman}
Joshua Lackman. \textit{Cohomology of Lie Groupoid Modules and the Generalized van Est Map.}
International Mathematics Research Notices, rnab027, (2021), https://doi.org/10.1093/imrn/rnab027
\bibitem{Lackman2}
Joshua Lackman. \textit{The van Est Map on Geometric Stacks.} arXiv:2205.02109
(2022).
\bibitem{Meinrenken} 
David Li-Bland and  Eckhard Meinrenken.
\textit{On the van Est homomorphism for Lie groupoids.} L’Enseignement Mathématique, , Vol. 61, (2014).
\bibitem{Masmoudi}
M. Masmoudi. \textit{Tangential formal deformations of the Poisson bracket and tangential
products on a regular Poisson manifold.} J. Geom. Phys. 9 (1992) 155-171.
\bibitem{Moyal}
 J. E. Moyal. \textit{Quantum mechanics as a statistical theory.} Mathematical Proceedings of the Cambridge Philosophical Society. Cambridge University Press (CUP). 45 (1) (1949): 99–124.
\bibitem{Murphy}
Gerard J. Murphy. \textit{C*-Algebras and Operator Theory.} Academic Press 1990.
\bibitem{zas}
A. Polishchuk and E. Zaslow. \textit{Categorical mirror symmetry: The Elliptic curve.}
Adv. Theor. Math. Phys. \textbf{2}, 443-470 (1998)
doi:10.4310/ATMP.1998.v2.n2.a9
[arXiv:math/9801119 [math.AG]].
\bibitem{rieffel}
M. Rieffel. \textit{Deformation quantization and operator algebras.} Proc. Sympos. Pure. Math., vol. 51, Amer. Math. Soc., Providence, RI 1990, pp. 411-423. 15.
\bibitem{rieffel2}
M. Rieffel. \textit{Lie group convolution algebras as deformation quantizations of linear Poisson
structures.} Amer. J. Math. 112 (1990), 657–685. MR1064995.
\bibitem{rieffel3}
M. Rieffel. \textit{Questions on quantization.} Operator Algebras and Operator Theory (1997, Shanghai), Contemp. Math. 228 (1998), 315–326, quant-ph/9712009.
\bibitem{Runde}
V. Runde. \textit{Amenable, Locally Compact Groups. In: Amenable Banach Algebras.} Springer Monographs in Mathematics. Springer, New York, NY.
\bibitem{vais}
Izu Vaisman. \textit{Geometric Quantization on Presymplectic Manifolds.} Monatshefte für Mathematik (1983)
Volume: 96, page 293-310
ISSN: 0026-9255; 1436-5081/e
\bibitem{Waldmann}
Stefan Waldmann. \textit{Convergence of products: From examples to a general framework.} EMS Surv. Math. Sci. 6 (2019), no. 1/2, pp. 1–31
\bibitem{wein}
A. Weinstein. \textit{Tangential deformation quantization and polarized
symplectic groupoids.} Deformation Theory and Symplectic
Geometry, S. Gutt, J. Rawnsley, and D. Sternheimer, eds.,
Mathematical Physics Studies 20, Kluwer, Dordrecht, 1997,
301-314.
\bibitem{weinstein}
Alan Weinstein. \textit{Symplectic Groupoids, Geometric Quantization, and Irrational Rotation Algebras.} In: Dazord, P., Weinstein, A. (eds) Symplectic Geometry, Groupoids, and Integrable Systems. Mathematical Sciences Research Institute Publications, vol 20,  (1991). Springer, New York, NY.
\bibitem{weinstein1}
Alan Weinstein and Ping Xu. \textit{Extensions of symplectic groupoids and quantization.}
Journal für die reine und angewandte Mathematik. Vol. 417, (1991) pp. 159-190.
\bibitem{Wigner}
E. P. Wigner. \textit{On the quantum correction for thermodynamic equilibrium.} Physical Review. 40 (5) (1932): 749–759.
\bibitem{willie}
Willie Wong. \textit{How to compute the asymptotics of this oscillatory integral?}  https://mathoverflow.net/q/439154.
\bibitem{zhuc}
Chenchang Zhu. \textit{Lie II theorem for Lie algebroids via higher Lie groupoids.}
Arxiv: Differential Geometry, (2010). arXiv:math/0701024v2 [math.DG]
\bibitem{Zachos}
Cosmas Zachos. \textit{Geometrical evaluation of star products}. J. Math. Phys. 41, 5129-5134 (2000) https://doi.org/10.1063/1.533395.
\end{thebibliography}
\end{document}